\documentclass[a4paper,12pt]{article}
\usepackage{amsmath}
\usepackage{amssymb,epic,eepic}
\usepackage[mathcal]{euscript}
\usepackage{epsfig}

\newcommand{\beq}{\begin{eqnarray*}}
\newcommand{\eeq}{\end{eqnarray*}}
\newcommand{\eps}{\epsilon}

\newcommand{\rz}{{\mathbb R}}
\newcommand{\gz}{{\mathbb Z}}

\newcommand{\sx}{{\mathbb S}}
\newcommand{\DD}{{\mathbb D}}
\newcommand{\KK}{{\mathbb K}}
\newcommand{\cn}{{\cal L}^{N+1}}

\newcommand{\cH}{{\cal H}}
\newcommand{\gh}{{\mathfrak h}}
\newcommand{\ee}{{\bf e}}

\newcommand{\OH}{{\cal O}}
\newcommand{\xx}{{x}}

\newcommand{\XX}{{X}}
\newcommand{\BX}{{\mathcal X}}
\newcommand{\PP}{{\mathcal P}}
\newcommand{\veps}{{\boldsymbol{\eps}}}
\newcommand{\rest}{\!\!\upharpoonright}

\newcommand{\BBox}{\hspace{1mm}\vrule height6pt width5.5pt depth0pt
                   \hspace{6pt}}

\newtheorem{prp}{Proposition}
\newtheorem{fct}{Fact}
\newtheorem{lem}{Lemma}
\newtheorem{thm}{Theorem}
\newtheorem{dfn}{Definition}
\newtheorem{gou}{Goursat Problem}

\title{Discrete and smooth orthogonal systems: $C^\infty$-approximation}
\author{A.I.\,Bobenko
 \thanks{Partially supported by the SFB 288 ``Differential Geometry and
 Quantum Physics'' and the DFG Research Center ``Mathematics for Key
 Technologies'' (FZT 86) in Berlin.}
\and D.\,Matthes
 \thanks{Supported by the SFB 288 ``Differential Geometry and
 Quantum Physics''.}
\and Yu.B.\,Suris$\;^\dag$}
\date{Institut f\"ur Mathematik, Technische Universit\"at Berlin,
 Str.~des 17.~Juni 136, 10623 Berlin, Germany.\\ \vspace{3mm}
 {\small E--Mails: {\tt bobenko}@{\tt math.tu-berlin.de},
 {\tt matthes}@{\tt math.tu-berlin.de},
 {\tt suris}@{\tt sfb288.math.tu-berlin.de}}}

\begin{document}
\maketitle


\section{Introduction and main results}
\label{sect:intro}

Triply orthogonal coordinate systems have attracted the attention
of mathematicians and physicists for almost two hundred years now.
Particular examples of them were already used by Leibniz and
Euler to evaluate multiple integrals in canonical coordinates, and
later Lam\'{e} and Jacobi carried out calculations in analytical
mechanics with the help of the famous elliptic coordinates.
The first general results on the geometry of triply orthogonal
systems date back to the 19th century -- like the famous theorem
of Dupin, saying that coordinate surfaces intersect along
curvature lines. The Lam\'{e} equations which analytically
describe the triply orthogonal systems were studied in detail by
Bianchi \cite{Bia} and Darboux \cite{Dar}.

Recently orthogonal systems came back into the focus of interest
in mathematical physics as an example of an integrable system.
Zakharov \cite{Zak} has shown how the Lam\'{e} equations can be
solved by the $\bar\partial$-method and constructed a variety of
explicit solutions with the help of the dressing method.
Algebro-geometric solutions of the Lam\'{e} equations were
constructed by Krichever \cite{Kri}. The recent interest to the
orthogonal coordinate systems is in particular motivated by their
applications to the theory of the associativity equations
\cite{Dub}.

The question of proper discretization of the classical models in
differential geometry became recently a subject of intensive study
(see, in particular, \cite{Bob,DoSa3,wolf}). Indeed one can
suggest various models in discrete geometry which have the same
continuous limit and nevertheless have quite different properties.
For a great variety of geometric problems described by integrable
equations it was found that {\em integrable discretizations} (i.e.
the discretizations preserving the integrability of the underlying
nonlinear system) preserve the characteristic geometric properties
of the problem. Moreover, it turns out that this discretization
can be described in pure geometric terms.

\begin{figure}[t]
\begin{center}
\epsfig{file=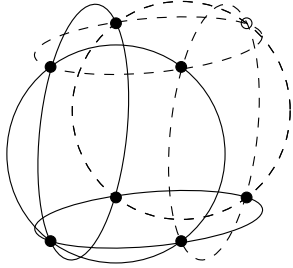,width=55mm} \caption{An elementary
hexahedron of a discrete orthogonal system} \label{fig:Kreise_a}
\end{center}
\end{figure}

\begin{figure}[ht]
\begin{minipage}[t]{180pt}
\epsfig{file=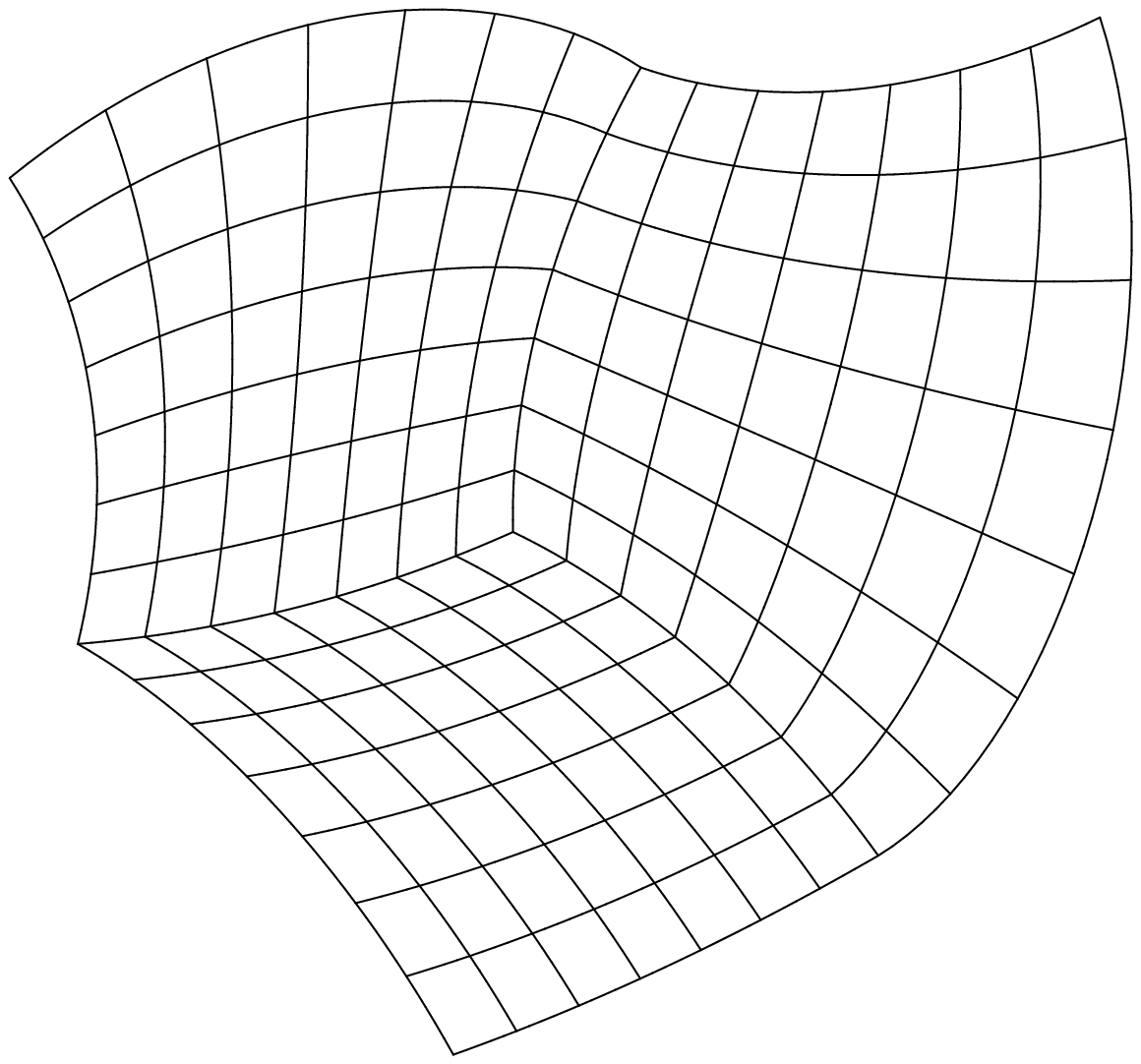,width=180pt}
\caption{Initial data for a smooth triply orthogonal system.}
\label{fig:Net}
\end{minipage}
\hfill
\begin{minipage}[t]{180pt}
\epsfig{file=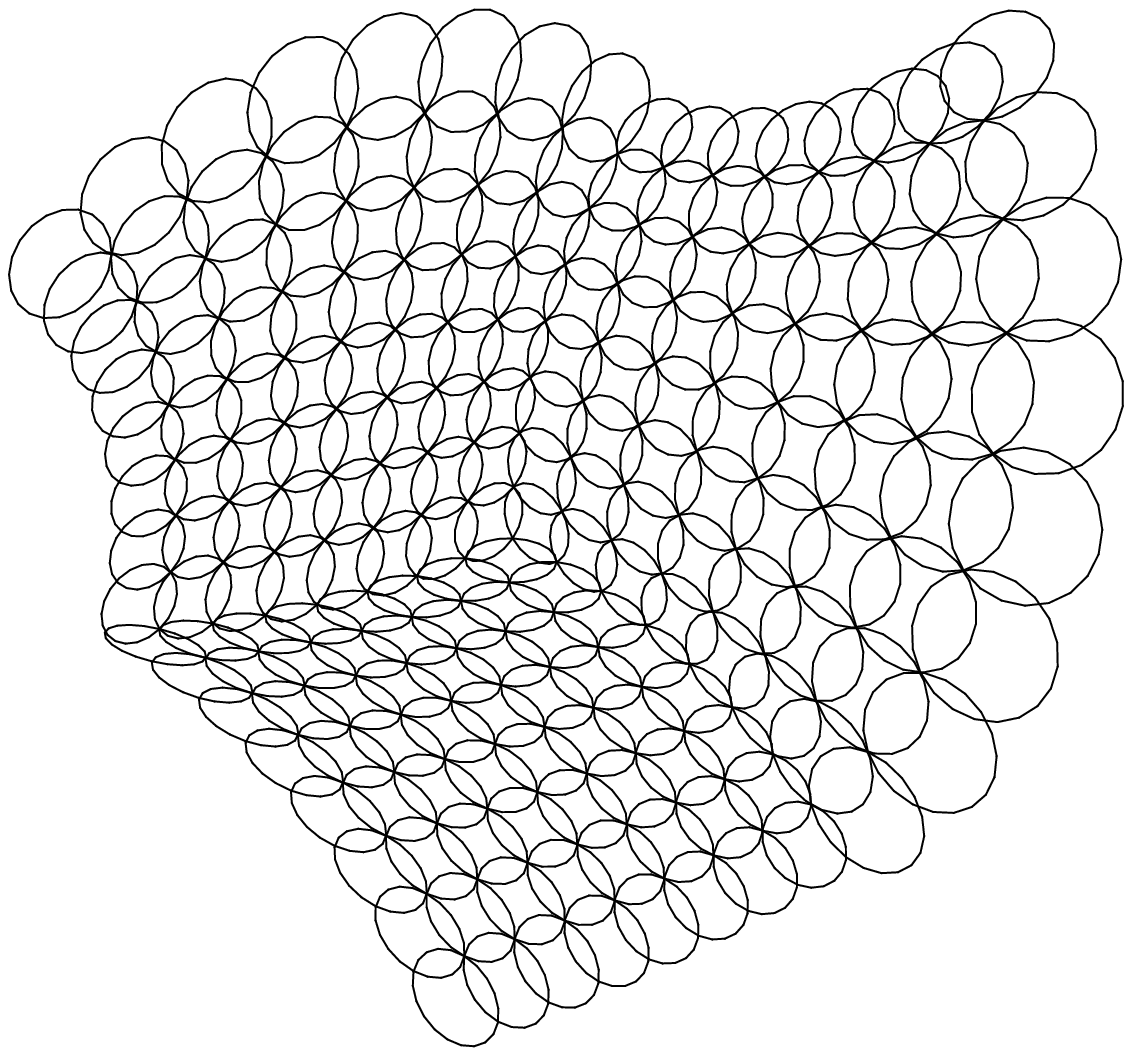,width=180pt}
\caption{Initial data for a discrete triply orthogonal system.}
\label{fig:Kreise_b}
\end{minipage}\hfill
\end{figure}

Such a discretization of the triply-orthogonal coordinate systems
was first suggested by one of the authors \cite{Bob96}. It is
based on the classical Dupin theorem. Since the circular nets are
known \cite{MPS,Nut} to correspond to the curvature line
parametrized surfaces, it is natural to define discrete triply
orthogonal systems as maps from ${\mathbb Z}^3$ (or a subset
thereof) to ${\mathbb R}^3$ with all elementary hexahedrons lying
on spheres. The neighboring spheres intersect along circles which
build the coordinate discrete curvature line nets (see
Fig.\ref{fig:Kreise_b}). Doliwa and Santini \cite{DoSa1} made the
next crucial step in the development of the theory. They considered
discrete orthogonal systems as a reduction of discrete conjugated systems
\cite{DoSa2}, generalized
them to arbitrary dimension and proved their geometric
integrability based on the classical Miguel theorem \cite{Ber}.
The latter claims that provided the seven (black) points lie on
circles as shown in Fig.\ref{fig:Kreise_a}, the three dashed
circles intersect in a common (white) point. This implies that a
discrete triply-orthogonal system is uniquely determined by its
three coordinate circular nets (as shown in
Fig.\ref{fig:Kreise_b}). This is a well-posed initial value
problem for triply orthogonal nets.

Later on discrete orthogonal systems were treated by analytic
methods of the theory of solitons. Algebro-geometric solutions
\cite{Volv1} as well as the $\bar\partial$-method \cite{DoSa4}
were discretized preserving all the
symmetries of the smooth case.

There is a common belief that the smooth theory can be obtained as
a limit of the corresponding discrete one, which can be treated
also in this case as a practical geometric numerical scheme for
computation of smooth surfaces. Moreover, this discrete master
theory has an important advantage: surfaces and their
transformations appear exactly in the same way as sublattices of
multidimensional lattices. This description immediately provides
the statements about the permutability of the corresponding
(B\"acklund-Darboux) transformations, -- an important property which
is nontrivial to see and sometimes even difficult to check in the
smooth theory.

Until recently there were no rigorous mathematical statements
supporting the observation about the classical limit of the
discrete theories. The first step in closing this gap was made in
our paper \cite{our} where general convergence results were proven
for a geometric numerical scheme for a class of nonlinear
hyperbolic equations. The first geometric example covered by this
theory was smooth and discrete surfaces with constant negative
Gaussian curvature described by the sine-Gordon equation.

In the present paper we formulate initial value problems and prove
the corresponding convergence results for conjugate and orthogonal
coordinate systems, which are described by three-dimensional
partial differential/difference equations. One should note here
that despite the known geometric discretization, a convenient
analytic description of discrete orthogonal systems converging to
the classical description of the smooth case was missing. In
particular this was the problem with discrete Lam\'{e} equations (cf.
\cite{Volv1}). The Clifford algebra description
of smooth and discrete orthogonal systems suggested in \cite{s&u} turned
out to be optimal for our purposes. This description provides us with the
proper discrete analogs of the continuous quantities and their
equations.

In particular we introduce the discrete Lam\'{e} system which is
the starting point for further analytical investigations. We write
these equations in the form of a discrete hyperbolic system and
formulate the corresponding Cauchy problem. As a by-product, we
rewrite the classical Lam\'{e} system in the hyperbolic form and
pose a Cauchy problem for classical orthogonal coordinate systems,
which seems to be new.

The central theorems of this paper are on the convergence of discrete
orthogonal coordinate systems to a continuous one. In particular
we prove the following statement.
\begin{thm}
  \label{thm:cDupin}
  Assume three umbilic free immersed surfaces ${\cal F}_1$, ${\cal F}_2$
  and ${\cal F}_3$ in $\rz^3$ intersect pairwise along their curvature lines,
  $\XX_1={\cal F}_2\cap{\cal F}_3$ etc. Then:
  \begin{enumerate}
  \item On a sufficiently small box $B=[0,r]^3\subset\rz^3$,
    there is an orthogonal coordinate system $\xx:B\rightarrow\rz^3$,
    which has ${\cal F}_1$, ${\cal F}_2$
    and ${\cal F}_3$ locally as images of the coordinate planes,
    with each curve $\XX_i$ being parametrized locally over
    the $i$-th axis by arc-length. The
    orthogonal coordinate system $\xx$ is uniquely determined by
    these properties.
  \item Denote by $B^\eps=B\cap(\eps\gz)^3$ the grid of mesh size $\eps$
    inside $B$. One can define a family $\{\xx^\eps:B^\eps\rightarrow\rz^3\}$
    of discrete orthogonal systems, parametrized by $\eps>0$, such that the
    approximation error decays as
    \begin{equation*}
      \sup_{\xi\in B^\eps}|\xx^\eps(\xi)-\xx(\xi)|<C \eps.
    \end{equation*}
    The convergence is $C^\infty$: all partial difference
    quotients of $\xx^\eps$ converge with the same rate to the respective
    partial derivatives of $\xx$.
  \end{enumerate}
\end{thm}

We also answer the question of {\it how to construct} the
approximating family $\{\xx^\eps\}$: One solves the discrete
Lam\'{e} system with discrete Cauchy data determined by ${\cal
F}_1$, ${\cal F}_2$ and ${\cal F}_3$. The solution $\xx^\eps$ is
then calculated in a process that consists of
$\OH(\eps^{-3})$-many evaluations of the expressions in the
discrete Lam\'{e} equations. The scheme can be -- and has been --
easily implemented on a computer, it needs a precision of $\eps^2$
for the calculations and is robust against $C^\infty$-deviations
of order $\eps$ in the initial data.

Actually, the results presented here go further:
\begin{itemize}
\item
We do not restrict ourselves to the case of three dimensions, but
allow arbitrary dimensions, possibly different for domain and
target space. In particular, we prove the corresponding approximation
results for curvature line parametrized surfaces.
\item
One can keep some of the directions of the orthogonal system
discrete while the others become continuous. This way, the
classical Ribaucour transformation is obtained very naturally. The
approximation result holds simultaneously for the system and its
Ribaucour transforms.
\item
The results about the permutability of the Ribaucour
transformations are elementary in the discrete setup. Our
approximation theorem implies the corresponding claim for
classical orthogonal systems (Theorem \ref{thm: rib permut}).
The first part of this theorem is due to Bianchi \cite{Bia0}, the second
part was proved in \cite{Tsa} using Bianchi's analytic description
of the Ribaucour transformations.
\item
In addition to orthogonal systems, we also treat, in a completely
analogous manner, conjugate nets. We pose a Cauchy problem for
them and show that any continuous conjugate net, along with its
transformations, can be $C^\infty$-approximated by a sequence of
discrete conjugate nets. The latter are quadrilateral lattices introduced
by Doliwa and Santini \cite{DoSa2}.
\end{itemize}

The paper is organized as follows.
We start in Sect. \ref{sect:general} with a short presentation of the
concepts developed in our recent paper \cite{our}, stating the main result
about an approximation of solutions of a hyperbolic system of partial
differential equations by solutions of the corresponding partial difference
system.
In Sect. \ref{sect conj nets} we apply this theorem to conjugate systems. We
formulate a Cauchy problem for smooth and discrete conjugate
systems, define the Jonas transformation and prove the
convergence.
Sect. \ref{sect otho systems} is dedicated to the description of orthogonal
systems in M\"obius geometry. Sect. \ref{sect: ortho approx}
contains the approximation results for orthogonal systems, a part of
which was summarized above.


\section{The approximation theorem}
\label{sect:general}
We start with a study of hyperbolic systems of partial difference
and partial differential equations. A general convergence theorem is
formulated, stating that solutions of difference equations approximate,
under certain conditions, solutions of differential equations, provided
the approximation of equations themselves takes place.

\subsection{Notations}

For a vector $\veps=(\eps_1,\ldots,\eps_M)$ of $M$ non-negative
numbers, we define the lattice
\begin{equation}\label{eq:lattice}
  \PP^\veps=(\eps_1\gz)\times\ldots\times(\eps_M\gz)\subset\rz^M.
\end{equation}
We refer to $\eps_i$ as the {\it mesh size} of the grid $\PP^\veps$ in the
$i$-th direction; a direction with a vanishing mesh size $\eps_i=0$
is understood as a {\it continuous} one, i.e. in this case the corresponding
factor $\eps_i\gz$ in (\ref{eq:lattice}) is replaced by $\rz$. For instance,
the choice $\veps=(0,\ldots,0)$ yields $\PP^\veps=\rz^M$. We use also the
following notation for $n$-dimensional sublattices of $\PP^\veps$ $(n\le M)$:
\begin{equation}\label{eq:sublattice}
  \PP^\veps_{i_1,\ldots,i_n}=\left\{\xi=(\xi_1,\ldots,\xi_M)\in\PP^\veps
  \mid \xi_i=0\mbox{ if }i\notin\{i_1,\ldots,i_n\}\right\}.
\end{equation}
Considering convergence problems, we will mainly deal with the following
two situations:
\begin{enumerate}
\item $\eps_i=\eps$ for all $i=1,\ldots,M$
  (all directions become continuous in the limit $\eps\to 0$);
\item $\eps_i=\eps$ for $i=1,\ldots,m$; $\eps_i=1$ for $i=m+1,\ldots,M$
  (the last $M-m$ directions are kept discrete in the limit).
\end{enumerate}
Thus, we always have one small parameter $\eps\geq 0$ only. In our geometric
problems, the first situation corresponds to the approximation of conjugate
and orthogonal systems themselves, while the second situation corresponds to
the approximation of such systems along with their transformations.
The underlying difference equations for discretizations of
conjugate (resp. orthogonal) systems, and for their transformations
will be the same, the only
distinction being in the way the continuous limit is performed. We will
distinguish between these two situations by setting $M=m+m'$, where $m'=0$
and in the first situation, while $m'>0$ in the second one.

Our considerations are local, therefore the encountered functions
are usually defined on compact domains, like the cube
\begin{equation}\label{eq:ball c}
B(r)=\left\{\xi\in\rz^M \mid 0\leq\xi_i\leq r\mbox{ for }i\leq M \right\}
\end{equation}
in the continuous context, and the following compact subsets
of $\PP^\veps$ in the lattice context. In the situation 1 ($m'=0$):
\begin{equation}\label{eq:ball 1}
  B^\eps(r)=\left\{\xi\in\PP^\veps \mid 0\leq\xi_i\leq r\mbox{ for }i\leq M
  \right\};
\end{equation}
in the situation 2 ($m'>0$):
\begin{equation}\label{eq:ball 2}
  B^\eps(r)=\left\{\xi\in\PP^\veps \mid 0\leq\xi_i\leq r
  \mbox{ for }1\le i\leq m,\;\;\xi_i\in\{0,1\}
  \mbox{ for }m<i\le M\right\}.
\end{equation}
We will often use the notation
\begin{equation}\label{eq:ball 0}
  B_0^\eps(r)=\left\{\xi\in(\eps\gz)^m \mid
  0\leq\xi_i\leq r\mbox{ for }1\le i\leq m \right\},
\end{equation}
so that in the latter case $B^\eps(r)=B_0^\eps(r)\times\{0,1\}^{m'}$.

For lattice functions $u:\PP^\veps\mapsto\BX$ with values in some Banach
space $\BX$, we use the shift operators  $\tau_i$ and the partial difference
operators $\delta_i=(\tau_i-{\bf 1})/\eps_i$ defined as
\begin{equation}\label{eq:shifts}
  (\tau_i u)(\xi)=u(\xi+\eps_i\ee_i),\qquad
  (\delta_i u)(\xi)=\frac{1}{\eps_i}
  \left(u(\xi+\eps_i\ee_i)-u(\xi)\right);
\end{equation}
here $\ee_i$ is the $i$th unit vector.
For $\eps_i=0$, define $\tau_i={\bf 1}$ and $\delta_i=\partial_i$.
Multiple shifts and higher order partial differences are written with the
help of multi-indices $\alpha=(\alpha_1,\ldots,\alpha_M)\in{\mathbb N}^M$:
\begin{equation}\label{eq:mult shifts}
  \tau^\alpha=\tau_1^{\alpha_1}\cdots\tau_M^{\alpha_M},\qquad
  \delta^\alpha=\delta_1^{\alpha_1}\cdots\delta_M^{\alpha_M},\qquad
  |\alpha|=\alpha_1+\ldots+\alpha_M.
\end{equation}
We use the following $C^\ell$-norms for functions $u:B^\eps(r)\rightarrow\BX$:
\begin{equation}\label{eq:C-norm}
  \|u\|_\ell=\sup\left\{|\delta^\alpha u(\xi)|_\BX\,\Big{|}\,
  \xi\in B^\eps(r-\eps|\alpha|);\,
  |\alpha|\leq\ell\right\};
\end{equation}
if $m'>0$, only the multi-indices $\alpha$ with $\alpha_i=0$ $\;(i>m)$
are taken into account in the last formula. Any function defined on
$B^\eps(r)$ with $\eps>0$
has finite $C^\ell$-norms of arbitrary order $\ell$; on the other hand, if
$u:B(r)\mapsto\BX$ is a function of a continuous argument which does not
belong to the class $C^\ell(B(r))$, then for its restriction to the lattice
points, $u^\eps=u\rest_{B^\eps(r)}$, the norm  $\|u^\eps\|_\ell$ will,
in general, diverge for $\eps\rightarrow 0$.
\begin{dfn}
 We say that a family of maps
 $\{u^\eps:B^\eps(r)\rightarrow\BX\}_{0<\eps<\eps_1}$ is
 $O(\eps)$-convergent in $C^\infty(B(r))$ to a smooth limit
 $u:B(r)\rightarrow\BX$, if there are finite real constants $K_\ell$ such that
 \begin{equation*}
    \|u^\eps-u\rest_{B^\eps(r)}\|_\ell\leq K_\ell\eps.
 \end{equation*}
 for any integer $\ell\geq 0$ and all positive $\eps<r/\ell$.
\end{dfn}

\subsection{Approximation for hyperbolic equations}
\label{sect hyperbolic}

Being interested in the convergence of solutions of difference equations
towards the solutions of differential ones, we restrict our attention to
two versions of the limit behaviour of the lattice $\PP^\veps$ described
above, i.e. $m'=0$ (all $\eps_i=\eps$) or $m'>0$ ($\eps_i=\eps$ for $1\le i
\le m$, and $\eps_i=1$ for $m+1\le i\le M=m+m'$).

Let $\BX=\BX_1\times\cdots\times\BX_K$ be a direct product of $K$
Banach spaces, equipped with the norm
\begin{equation}\label{eq:knorm}
  |u(\xi)|_{\BX}=\max_{1\le k\le K} |u_k(\xi)|_{{\BX}_k}.
\end{equation}
Consider a hyperbolic system of first-order partial difference ($\eps>0$) or
partial differential ($\eps=0$) equations for functions $u:\PP^\veps\mapsto\BX$:
\begin{equation} \label{eq:evolution}
  \delta_j u_k=f^\eps_{k,j}(u),
  \qquad 1\le k\le K,\quad j\in{\cal E}(k).
\end{equation}
Thus, for each $k=1,\ldots,K$, there is defined a subset of indices
${\cal E}(k)\subset\{1,\ldots,M\}$, called the set of {\it evolution
directions} for the component $u_k$; the complementary subset
${\cal S}(k)=\{1,\ldots M\}\setminus{\cal E}(k)$ is called the set
of {\it static directions} for the component $u_k$. The dependent variables
$u_k(\xi)=u_k(\xi_1,\ldots,\xi_M)$ are thought of as attached to the cubes of
dimension $\#({\cal S}_k)$ adjacent to the point $\xi$.

The functions $f_{k,j}^\eps$ are defined on open domains
${\DD}_{k,j}(\eps)\subset\BX$ and are of class $C^\infty$ there.
We suppose that the following conditions are satisfied.
\begin{itemize}
\item[(F)] The domains ${\DD}_{k,j}(\eps)$ for $1\le j\le m$
 (that is for all directions which become continuous in the limit $\eps\to 0$)
 exhaust $\BX$ as $\eps\to 0$, i.e. any compact subset $\KK\subset\BX$
 is contained in $\DD_{k,j}(\eps)$ for all $k$, for all
 $j\in{\cal E}(k)$, $j\leq m$, and for all $\eps$ small enough.
 Similarly, there exists an open set $\DD\subset\BX$ such that
 the domains ${\DD}_{k,j}(\eps)$ for $m+1\le j\le M$ corresponding
 to the directions remaining discrete in the limit,
 exhaust $\DD$ as $\eps\to 0$.
 The functions $f^\eps_{k,j}$ locally converge to the respective
 $f^0_{k,j}$ in $C^\ell(\BX)$, resp. in $C^\ell(\DD)$, for any $\ell$,
 with the rate of convergence $O(\eps)$.
\end{itemize}
The existence of the subset $\DD$ reflects the fact that transformations of
conjugate, resp. orthogonal, systems are not defined for all initial data.
The set $\DD$ consists of ``good'' initial values, where the transformations
are well defined.

\begin{gou}
  \label{gou:abstract}
  Given functions $U^\eps_k$ on $\PP^\eps_{{\cal S}(k)}$
  (the subspace spanned by the static directions of $u_k$),
  and $r>0$, find a function $u^\eps:B^\eps(r)\rightarrow\BX$ that
  solves the system (\ref{eq:evolution}) and satisfies the initial
  conditions $u^\eps_k=U^\eps_k$ on $\PP^\eps_{{\cal S}(k)}\cap B^\eps(r)$.
\end{gou}

In order for this problem to admit a unique solution for $\eps>0$, this has
to be the case
just for one elementary $M$-dimensional cube $B^\eps(\eps)$. The corresponding
condition is called {\em consistency}, and it turns out to be necessary and
sufficient also for the global (formal) solvability of the above Goursat
problem. The consistency condition can be expressed as
\begin{equation}\label{comp}
\delta_{j}(\delta_{i}u_k)=\delta_{i}(\delta_{j}u_k),
\end{equation}
for any choice of $i\neq j$ from the respective ${\cal E}_k$, if one
writes out all the involved $\delta u$'s according to the equations
(\ref{eq:evolution}). The first step in writing this out reads:
\begin{equation}\label{comp1}
\delta_{j}(f_{k,i}(u))=\delta_{i}(f_{k,j}(u)),
\end{equation}
In order for the left--hand side (say) of this
equation to be well--defined, the function $f_{k,i}$ is allowed to depend
only on those $u_{\ell}$ for which $j\in{\cal E}_{\ell}$. This has to hold
for all $j\in{\cal E}_k$, $j\neq i$, so we come to the condition that
$f_{k,i}$ depends only on those $u_{\ell}$ for which
${\cal E}_k\setminus\{i\}\subset{\cal E}_\ell$. If this is satisfied,then the
final form of the consistency condition reads:
\begin{equation}\label{comp2}
\eps_j^{-1}\Big(f_{k,i}(u+\eps_jf_j(u))-f_{k,i}(u)\Big)=
\eps_i^{-1}\Big(f_{k,j}(u+\eps_if_i(u))-f_{k,j}(u)\Big),
\end{equation}
where we use the abbreviation $f_i(u)$ for the $K$--vector whose $k$th
component is $f_{k,i}(u)$, if $i\in{\cal E}_k$, and is not defined otherwise.

It is readily seen that, if one prescribes arbitrarily
the values of all variables $u_k$ at $\xi=0$ (recall, $u_k(0)$ is actually
associated with the respective static facet of dimension
$\#({\cal S}_k)$ adjacent to the point $\xi=0$), then the system
(\ref{eq:evolution}) defines uniquely the values for all variables on
the remaining facets of the unit cube if and only if it is consistent.
The following proposition shows that the local consistency assures also
the global solvability of the Goursat problem; the only obstruction
could appear if the solution would hit the set where the functions $f_{k,j}$
are not defined.

\begin{prp}\label{lem:dcn}
  Let $\eps>0$, and assume that the system (\ref{eq:evolution}) is consistent.
  Then either the Goursat problem \ref{gou:abstract} has a unique solution
  $u^\eps$ on the whole of $B^\eps(r)$, or there is a subdomain
  $\Omega\subset B^\eps(r)$ such that the Goursat problem \ref{gou:abstract}
  has a unique solution $u^\eps$ on $\Omega$, and
  $u^\eps(\xi)\notin{\DD}(\eps)$ at some point $\xi\in\Omega$.
\end{prp}

Finally, we turn to the whole family ($0\le \eps\leq \eps_0$) of Goursat
problems for equations (\ref{eq:evolution}) with the respective data
$U^\eps_k$. The following result has been proven in our recent paper
\cite{our}.
\begin{thm}
  \label{thm:approx}
  Let a family of Goursat problems \ref{gou:abstract} be given.
  Suppose that they are consistent for all $\eps>0$, and satisfy
  the condition (F).
  \begin{enumerate}
  \item
    Let the functions $U^0_1,\ldots,U_K^0$ be smooth,
    and, if $m'>0$, assume that the point $(U^0_1(0),\ldots,U^0_K(0))$
    belongs to the subset $\DD$. Then there exists $r>0$ such that the
    Goursat problem \ref{gou:abstract} at $\eps=0$ has a unique solution
    $u:B(r)\rightarrow\BX$.
  \item
    If the Goursat data $U^\eps_k$ locally $O(\eps)$-converge
    in $C^\infty$ to the functions $U_k^0$, then there exists
    $\eps_1\in(0,\eps_0)$ such that the Goursat problems \ref{gou:abstract}
    are solvable on $B^\eps(r)$ for all $0<\eps<\eps_1$, and the solutions
    $u^\eps:B^\eps(r)\rightarrow\BX$ converge to $u:B(r)\rightarrow\BX$
    with the rate $O(\eps)$ in $C^\infty(B(r))$.
  \end{enumerate}
\end{thm}


\section{Conjugate nets}
\label{sect conj nets}

We give definitions of continuous and discrete conjugate nets and their
transformations, formulate corresponding Goursat problems, and prove
convergence of discrete conjugate nets towards continuous ones.

\subsection{Basic definitions}

Recall that $M=m+m'$. In Definition \ref{dfn:cn} below it
is supposed that $m'=0$, $M=m$, while in Definition \ref{dfn:jp} it
is supposed that $m'=1$, $M=m+1$.
\begin{dfn}\label{dfn:cn}
  A map $\xx:B_0(r)\mapsto\rz^N$ is called an $m$-dimensional
  {\em conjugate net} in $\rz^N$, if $\partial_i\partial_j\xx\in{\rm span}
  (\partial_i\xx,\partial_j\xx)$ at any point $\xi\in B_0(r)$ for all pairs
  $1\leq i\neq j\leq m$, i.e. if there exist functions
  $c_{ij},c_{ji}:B_0(r)\mapsto\rz$ such that
  \begin{equation} \label{eq:ccnproperty}
    \partial_i\partial_j\xx=c_{ji}\partial_i\xx+c_{ij}\partial_j\xx\,.
  \end{equation}
\end{dfn}
\begin{dfn}\label{dfn:jp}
  A pair of $m$-dimensional conjugate nets
  $\xx,\xx^+:B_0(r)\mapsto\rz^N$ is called a {\em Jonas pair}, if
  three vectors $\partial_i\xx$, $\partial_i\xx^+$ and $\delta_Mx=\xx^+-\xx$
  are coplanar at any point $\xi\in B_0(r)$ for all $1\le i\le m$, i.e.
  if there exist functions $c_{iM}, c_{Mi}:B_0(r)\rightarrow\rz$ such that
  \begin{equation} \label{eq:cjtproperty}
  \partial_i\xx^+-\partial_i\xx=c_{Mi}\partial_i\xx+c_{iM}(\xx^+-\xx)\,.
  \end{equation}
\end{dfn}
{\bf Remarks.}
\begin{itemize}
\item $\xx^+$ is also called a {\it Jonas transformation} of $\xx$.
One can iterate these transformations and obtain sequences
$(\xx^{(0)},\xx^{(1)},\ldots,\xx^{(T)})$ of conjugate nets
$\xx^{(t)}:B_0(r)\rightarrow\rz^N$.
It is natural to think of such sequences as of conjugate systems with $m$
continuous and one discrete direction. We will see immediately that the
notion of fully discrete conjugate systems puts all directions on an
equal footing.
\item A {\it Combescure transformation} $\xx^+$ of $\xx$ is
a Jonas transformation for which vectors $\partial_i\xx^+$ and
$\partial_i\xx$ are parallel, $i=1,\ldots,m$. This class of transformations
is singled out by requiring $c_{iM}=0$ in equation (\ref{eq:cjtproperty}).
\end{itemize}

Discrete conjugate nets were introduced by Doliwa and Santini \cite{DoSa2}.
\begin{dfn}\label{dfn:dcn}
  A map $\xx:B^\eps(r)\mapsto\rz^N$ is called an
  $M$-dimensional {\em discrete conjugate net} in $\rz^N$, if the four
  points $\xx$, $\tau_i\xx$, $\tau_j\xx$, and $\tau_i\tau_j\xx$ are coplanar
  at any $\xi\in B^\eps(r)$ for all pairs $1\leq i\neq j\leq M$, i.e. if
  there exist functions $c_{ij}:B^\eps(r)\rightarrow\rz$ such that
  \begin{equation} \label{eq:dcnproperty}
    \delta_i\delta_j\xx=c_{ji}\delta_i\xx+c_{ij}\delta_j\xx.
  \end{equation}
  If $m'=1$, $M=m+1$, then the $M$-dimensional discrete conjugate net,
  considered as a pair of functions
  $\xx,\xx^+:B_0^\eps(r)\mapsto\rz^N$, is also
  called a {\em Jonas pair} of $m$-dimensional discrete conjugate nets.
\end{dfn}

\subsection{Hyperbolic equations for conjugate nets}

Introducing $M$ new functions $w_i:B^\eps(r)\rightarrow\rz^N$, we can
rewrite (\ref{eq:dcnproperty}) as system of first order:
\begin{eqnarray}
  \delta_i\xx & = & w_i\,, \label{eq:dcn1}\\
  \delta_i w_j & = & c_{ji}w_i+c_{ij}w_j\,,  \quad i\neq j, \label{eq:dcn2}\\
  \delta_i c_{kj} & = &
  (\tau_j c_{ik})c_{kj}+(\tau_j c_{ki})c_{ij}-(\tau_i c_{kj})c_{ij}, \quad
   i\neq j\neq k\neq i. \quad\label{eq:dcn3}
\end{eqnarray}
For a given discrete conjugate net $\xx:B^\eps(r)\rightarrow\rz^N$,
eq. (\ref{eq:dcn1}) defines the functions $w_j$,
then eq. (\ref{eq:dcn2}) reflects the property (\ref{eq:dcnproperty}),
and eq. (\ref{eq:dcn3}) is just a transcription of the compatibility
condition $\delta_i(\delta_j w_k)=\delta_j(\delta_i w_k)$. Conversely,
for any solution of (\ref{eq:dcn1})--(\ref{eq:dcn3}) the map
$\xx:B^\eps(r)\rightarrow\rz^N$ has the defining property
(\ref{eq:dcnproperty}) and thus is a discrete conjugate net. So,
discrete conjugate nets are in a one-to-one correspondence to
solutions of the system (\ref{eq:dcn1})-(\ref{eq:dcn3}).

This system almost suits the framework of Sect. \ref{sect hyperbolic};
the only obstruction is the implicit nature of the equations
(\ref{eq:dcn3}) (their right--hand sides depend on the shifted variables
like $\tau_ic_{kj}$ which is not allowed in (\ref{eq:evolution}).
We return to this point later, and for a moment we handle the system
(\ref{eq:dcn1})-(\ref{eq:dcn3}) just as if it would belong to
the class (\ref{eq:evolution}). In this context, we have to assign to every
variable a Banach space and static/evolution directions.
Abbreviating $w=(w_1,\ldots,w_M)$ and $c=(c_{ij})_{i\neq j}$, we set
\begin{equation*}
  \BX=\rz^N\{\xx\}\times(\rz^N)^M\{w\}\times\rz^{M(M-1)}\{c\},
\end{equation*}
with the norm from (\ref{eq:knorm}). We assign to $\xx$ no static directions,
to $w_i$ the only static direction $i$, and to $c_{ij}$ two static directions
$i$ and $j$. The corresponding Goursat problem is now formulated as follows.
\begin{gou}\label{gou:dcn}
{\bf (for discrete conjugate nets).}
  Given:
  \begin{itemize}
   \item a point $\XX^\eps\in\rz^N$,
   \item $M$ functions $W^\eps_i:\PP^\eps_i\mapsto\rz^N$ on the respective
      coordinate axes $\PP^\eps_i$, and
   \item $M(M-1)$ functions $C^\eps_{ij}:\PP^\eps_{ij}\mapsto\rz$ on the
     coordinate planes $\PP^\eps_{ij}$,
  \end{itemize}
  find a solution $(\xx^\eps,w^\eps,c^\eps):B^\eps(r)\mapsto\BX$
  to the equations (\ref{eq:dcn1})-(\ref{eq:dcn3}) satisfying the initial
  conditions
  \begin{equation}
   \xx^\eps(0)=\XX^\eps,\qquad
    w_i^\eps\rest_{\PP_i^\eps}=W_i^\eps, \qquad
    c^\eps_{ij}\rest_{\PP_{ij}^\eps}=C_{ij}^\eps.
   \end{equation}
\end{gou}

We discuss now the consistency of the discrete hyperbolic system
(\ref{eq:dcn1})--(\ref{eq:dcn3}). First of all, note that for any triple of
pairwise different indices $(i,j,k)$ the equations of this system involving
these indices only form a closed subset. In this sense, the system
(\ref{eq:dcn1})--(\ref{eq:dcn3}) consists of three-dimensional building blocks.
One says also that it is {\it essentially three-dimensional}. Correspondingly,
in the case $M=3$ consistency is not an issue: all conditions to be verified
are automatically taken into account by the very construction of the
equations (\ref{eq:dcn1})--(\ref{eq:dcn3}).

In this case $M=3$ the initial data for an elementary cube are:
a point $\xx(0)$,
three vectors $w_1(0)$, $w_2(0)$, $w_3(0)$ which are thought of as attached
to the edges incident to the point $\xi=0$ and parallel to the correspondent
axes, and six numbers $c_{ij}(0)$, $1\le i\neq j\le 3$, which are thought of
as attached to the plaquettes incident to the point $\xi=0$ and parallel to
the coordinate planes $\PP_{ij}$. From these data, one constructs first with
the help of (\ref{eq:dcn1}) the three points $\tau_i\xx=\xx(0)+\eps_iw_i(0)$,
then one calculates $\tau_i w_j$ for $i\neq j$ with the help of
(\ref{eq:dcn2}), and further $\tau_j\tau_i\xx$ by means of (\ref{eq:dcn1}).
Notice that the consistency $\tau_j\tau_i\xx=\tau_i\tau_j\xx$ is assured
by the symmetry of the right--hand side of (\ref{eq:dcn2}) with respect to
the indices $i$ and $j$. Thus, the
first two equations (\ref{eq:dcn1}), (\ref{eq:dcn2}) allow us to determine
the values of $\xx$ in seven vertices of an elementary cube. The remaining
8-th one, $\tau_1\tau_2\tau_3\xx$, is now determined as follows. Calculate
$\tau_ic_{kj}$ from the (linearly implicit) equations (\ref{eq:dcn3}),
then calculate $\tau_j\tau_iw_k$ from (\ref{eq:dcn2}), and then
$\tau_1\tau_2\tau_3x$ from (\ref{eq:dcn1}). Again, the consistency requirement
$\tau_j\tau_iw_k=\tau_i\tau_jw_k$ is guaranteed, per construction, by
eqs. (\ref{eq:dcn3}). Geometrically, the need to solve the linearly implicit
equations for $\tau_ic_{kj}$ is interpreted as follows. The point
$\tau_1\tau_2\tau_3x$ is determined by the conditions that
it lies in three planes $\tau_i\Pi_{jk}$, $1\le i\le 3$, where
$\tau_i\Pi_{jk}$ is defined as the plane passing through three points
$(\tau_i\xx,\tau_j\tau_i\xx,\tau_k\tau_i\xx)$, or, equivalently, passing
through $\tau_i\xx$ and spanned by $\tau_iw_j$, $\tau_iw_k$. All three planes
$\tau_1\Pi_{23}$, $\tau_2\Pi_{13}$ and $\tau_3\Pi_{12}$ belong to the
three--dimensional space through $\xx(0)$ spanned by $w_i(0)$, $1\le i\le 3$,
and therefore generically they intersect at exactly one point. The above
mentioned linearly implicit system encodes just finding the intersection
point of three planes.

Turning now to the case $M>3$, we face a non--trivial consistency problem.
Indeed, the point $\tau_i\tau_j\tau_k\tau_{\ell}\xx$ can be constructed
as the intersection point of three planes $\tau_i\tau_j\Pi_{k\ell}$,
$\tau_i\tau_k\Pi_{j\ell}$, $\tau_i\tau_\ell\Pi_{jk}$, where the plane
$\tau_i\tau_j\Pi_{k\ell}$ is
defined as the plane passing through three points
$(\tau_i\tau_j\xx,\tau_i\tau_j\tau_k\xx,\tau_i\tau_j\tau_\ell\xx)$, or,
equivalently, passing through $\tau_i\tau_j\xx$ and spanned by
$\tau_i\tau_jw_k$, $\tau_i\tau_jw_\ell$. So, there are four alternative
ways to determine the point $\tau_i\tau_j\tau_k\tau_\ell\xx$, depending
on which index plays the role of $i$.
\begin{prp}\label{prop:dcn consistency}
The system (\ref{eq:dcn1})-(\ref{eq:dcn3}) is consistent for $M>3$.
\end{prp}
{\bf Proof.} The geometric meaning of this statement is that the above
mentioned four ways to obtain $\tau_i\tau_j\tau_k\tau_{\ell}\xx$ lead to
one and the same result. This can be proved by symbolic manipulations
with the equations (\ref{eq:dcn1})-(\ref{eq:dcn3}), but a geometric
proof first given in \cite{DoSa2} is much more transparent.
In the construction above, it is easy to understand that the plane
$\tau_i\tau_j\Pi_{k\ell}$ is the intersection of two three--dimensional
subspaces $\tau_i\Pi_{jk\ell}$ and $\tau_j\Pi_{ik\ell}$ in the
four--dimensional space through $\xx(0)$ spanned by $w_i(0)$, $w_j(0)$,
$w_k(0)$, $w_\ell(0)$, where the subspace $\tau_i\Pi_{jk\ell}$ is the one
through the four points
$(\tau_i\xx,\tau_i\tau_j\xx,\tau_i\tau_k\xx,\tau_i\tau_\ell\xx)$, or,
equivalently, the one through $\tau_i\xx$ spanned by $\tau_iw_j$, $\tau_iw_k$
and $\tau_iw_\ell$. Now the point $\tau_i\tau_j\tau_k\tau_{\ell}\xx$
can be alternatively described as the unique intersection point of the
four three--dimensional subspaces $\tau_i\Pi_{jk\ell}$, $\tau_j\Pi_{ik\ell}$,
$\tau_k\Pi_{ij\ell}$ and $\tau_{\ell}\Pi_{ijk}$ of one and the same
four--dimensional space. \BBox
\bigskip

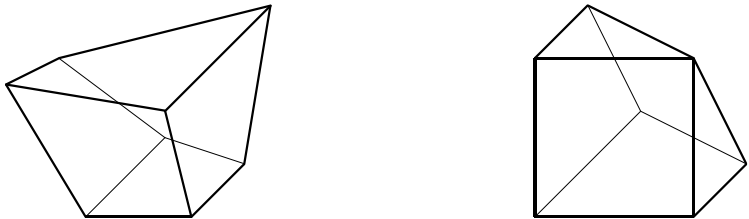
\begin{figure}[htbp]
\begin{picture}(300,100)(-50,0)
  \thicklines
  \put(30,0){\line(1,0){40}}
  \put(30,0){\line(-3,5){30}}
  \put(70,0){\line(1,1){20}}
  \put(70,0){\line(-1,4){10}}
  \put(0,50){\line(6,-1){60}}
  \put(0,50){\line(2,1){20}}
  \put(60,40){\line(1,1){40}}
  \put(20,60){\line(4,1){80}}
  \put(90,20){\line(1,6){10}}
  \thinlines
  \put(60,30){\line(3,-1){30}}
  \put(60,30){\line(-4,3){40}}
  \put(60,30){\line(-1,-1){30}}

  \thicklines
  \put(200,0){\line(1,0){60}}
  \put(200,0){\line(0,1){60}}
  \put(200,60){\line(1,0){60}}
  \put(260,0){\line(0,1){60}}
  \put(200,60){\line(1,1){20}}
  \put(260,0){\line(1,1){20}}
  \put(260,60){\line(-2,1){40}}
  \put(260,60){\line(1,-2){20}}
  \thinlines
  \put(200,0){\line(1,1){40}}
  \put(240,40){\line(-1,2){20}}
  \put(240,40){\line(2,-1){40}}
\end{picture}
\caption{A non-degenerate and a degenerate hexahedron.}
\label{fig:crashedpoly}
\end{figure}

Now, we turn to the above mentioned feature of the system
(\ref{eq:dcn1})--(\ref{eq:dcn3}), namely to its implicit nature.
The geometric background of this feature is the following. Though there is,
in general, a unique way to construct the 8-th vertex $\tau_i\tau_j\tau_kx$
of an elementary hexahedron out of the the seven vertices $x$, $\tau_ix$
and $\tau_i\tau_jx$, or, equivalently, out of $x$, $w_i$, $c_{ij}$, it
can happen that the resulting hexahedron is degenerate. This happens, e.g.,
if the point $\tau_i\tau_jx$ lies in the plane $\tau_k\Pi_{ij}$ through
$(\tau_kx,\tau_i\tau_kx,\tau_j\tau_kx)$. This is illustrated on Fig.
\ref{fig:crashedpoly}.
While the example on the left of Fig. \ref{fig:crashedpoly} is fine,
the right one does not have the combinatorics of a 3-cube
(the top-right edge has degenerated to a point). As a result, some of
the quantities $\tau_i c_{kj}$ in equation (\ref{eq:dcn3}) are not
well--defined; this is reflected by the fact that the equations
(\ref{eq:dcn3}) are implicit and have to be solved for $\tau_i c_{kj}$,
which is possible not for all values of $c$ but rather for those belonging
to a certain open subset.

To investigate this point, rewrite (\ref{eq:dcn3}) as
\begin{equation}
  \delta_i c_{kj}=
  F_{(ijk)}(c)+\eps_j(\delta_j c_{ik})c_{kj}+\eps_j(\delta_j c_{ki})c_{ij}
  -\eps_i(\delta_i c_{kj})c_{ij}
\end{equation}
Introducing vectors $\delta c$ and $F(c)$ with $M(M-1)(M-2)$ components
labelled by triples $(ijk)$ of pairwise distinct numbers $1\le i,j,k\le M$,
\begin{eqnarray*}
  \delta c_{(ijk)}&=&\delta_i c_{kj}\\
  F_{(ijk)}(c) & = & c_{ik}c_{kj}+c_{ki}c_{ij}-c_{kj}c_{ij},
\end{eqnarray*}
we present the above equations as
\begin{equation}
\label{eq:Qmatrix}
  ({\bf 1}-Q^\eps(c))\delta c=F(c),
\end{equation}
with a suitable matrix $Q^\eps(c)$. We need to find conditions for
${\bf 1}-Q^\eps(c)$ to be invertible.

\subsubsection{One conjugate net}
If $m'=0$, then all entries of $Q^\eps(c)$ are either zero or of the
form $\eps c_{ij}$. Therefore, on any compact ${\mathbb K}\subset\BX$ we
have: $\|Q^\eps(c)\|<C\eps$ with some $C=C({\mathbb K})>0$, hence
${\bf 1}-Q^\eps(c)$  is invertible for $\eps<\eps_0=1/C$, and the
the inverse matrix $({\bf 1}-Q^\eps(c))^{-1}$ is a smooth function
on ${\mathbb K}$ and $O(\eps)$-convergent in $C^\infty({\mathbb K})$
to ${\bf 1}$. So, the system (\ref{eq:Qmatrix}) is solvable, and
\begin{equation*}
  \delta_i c_{kj}=F_{(ijk)}(c)+O(\eps)=
  c_{ik}c_{kj}+c_{ki}c_{ij}-c_{kj}c_{ij}+O(\eps),
\end{equation*}
where the constant in $O(\eps)$ is uniform on any compact set.

The limiting equations for (\ref{eq:dcn1})-(\ref{eq:dcn3}) in this case are:
\begin{eqnarray}
   \partial_i\xx & = & w_i   \label{eq:ccn1.}\\
   \partial_i w_j & = & c_{ji}w_i+c_{ij}w_j   \label{eq:ccn2.}\\
   \partial_i c_{kj} & = & c_{ik}c_{kj}+c_{ki}c_{ij}-c_{kj}c_{ij}.
   \label{eq:ccn3.}
\end{eqnarray}
For any solution $(\xx,w,c):B_0(r)\rightarrow\BX$ of this system
the function $\xx:B_0(r)\rightarrow\rz^N$ is a conjugate net,
because equations (\ref{eq:ccn1.}) and (\ref{eq:ccn2.}) yield
the defining property (\ref{eq:ccnproperty}). Conversely, if
$\xx:B_0(r)\rightarrow\rz^N$ is a conjugate net,
then equation (\ref{eq:ccn1.}) defines the vectors $w_i$,
these fulfill (\ref{eq:ccn2.}) since this is the defining property of
a conjugate net, and (\ref{eq:ccn3.}) results from equating
$\partial_i(\partial_j w_k)=\partial_j(\partial_i w_k)$. So, the following is
demonstrated:
\begin{lem}\label{lemma conj}
  If $m'=0$, then the Goursat problem \ref{gou:dcn} for discrete conjugate
  nets fulfills the condition (F). The limiting system  for $\eps=0$,
  consisting of (\ref{eq:ccn1.})-(\ref{eq:ccn3.}),
  describes continuous conjugate nets.
\end{lem}

\subsubsection{Jonas pair of conjugate nets}

In the case $m'=1$ the matrix ${\bf 1}-Q^\eps(c)$ is no longer a small
perturbation of the identity, since $\eps_M=1$. More precisely, the matrix
$Q^\eps(c)$ is block--diagonal, its $M \choose 3$ diagonal blocks
$Q^\eps_{\{ijk\}}(c)$ of size $6\times 6$ correspond to non--ordered
triples of pairwise distinct indices $\{i,j,k\}$, with rows and
columns labeled by the six possible permutations of these indices.
Blocks for which $i,j,k\neq M$ are of the form
considered before -- their entries are either zero or of the form
$\eps c_{ij}$. Blocks where, say, $k=M$, admit a decomposition
$Q^\eps_{\{ijM\}}(c)=A(c)+\eps B(c)$, where the matrices $A$ and $B$
do not depend on $\eps$. It is not difficult to calculate that
$\det({\bf 1}-A(c))=(1+c_{Mi})(1+c_{Mj})$.
So ${\bf 1}-Q^\eps_{\{ijM\}}(c)$ is invertible if $c_{Mi},c_{Mj}\neq -1$,
provided $\eps$ is small enough. The inverse is smooth in a neighborhood
any such point $c$, and $O(\eps)$-convergent to $({\bf 1}-A(c))^{-1}$.

As formal limit of (\ref{eq:dcn1}) and (\ref{eq:dcn2}),
we obtain $(i\neq j)$:
\begin{eqnarray}
    \partial_i\xx & = & w_i \qquad (1\leq i\leq m),  \label{eq:cnlimit1}\\
    \delta_M\xx & = & w_M, \label{eq:cnlimit2}\\
    \partial_i w_j & = & c_{ji}w_i+c_{ij}w_j \qquad (1\leq i\leq m),
    \label{eq:cnlimit3}\\
    \delta_M w_j & = & c_{jM} w_M+c_{Mj} w_j. \label{eq:cnlimit4}
\end{eqnarray}
There are also four different limits of equation (\ref{eq:dcn3}),
depending on which directions $i$ and $j$ are kept discrete.
We shall not write them down; they are easily reconstructed as the
compatibility conditions for (\ref{eq:cnlimit3}), (\ref{eq:cnlimit4}),
like $\delta_M(\partial_i w_k)=\partial_i(\delta_M w_k)$ etc.
Comparing (\ref{eq:cnlimit1})--(\ref{eq:cnlimit4}) with the definition of
the Jonas transformation (\ref{eq:cjtproperty}), we see that the following is
demonstrated.

\begin{lem}\label{lemma jonas}
If $m'=1$, then the Goursat problem \ref{gou:dcn} for discrete conjugate nets
also fulfills  the condition (F), with the domain
  \begin{equation*}
    \DD=\{(\xx,w,c)\in\BX\mid c_{Mi}\neq -1\;\mbox{ for }\;1\leq i\leq m\}.
  \end{equation*}
  The limiting system at $\eps=0$ describes Jonas pairs
  of continuous conjugate nets.
\end{lem}

\subsection{Approximation theorems for conjugate nets}

We are now ready to formulate and prove the main results of this chapter.
\begin{thm} \label{thm:cnapprox}
{\bf (approximation of a conjugate net).}
  Let there be given:
  \begin{itemize}
  \item $m$ smooth curves $\XX_i:\PP_i\mapsto\rz^N$, $i=1,\ldots,m$,
  intersecting at a common point $\XX=\XX_1(0)=\cdots=\XX_m(0)$;
  \item for each pair $1\leq i<j\leq m$, two smooth functions
  $C_{ij},C_{ji}:\PP_{ij}\mapsto\rz$.
  \end{itemize}
  Then, for some $r>0$:
  \begin{enumerate}
  \item   There is a unique conjugate net $\xx:B_0(r)\rightarrow\rz^N$
  such that
  \begin{equation}\label{cn:goursat}
  x\rest_{\PP_i}=X_i\,,\quad c_{ij}\rest_{\PP_{ij}}=C_{ij}\,.
  \end{equation}
  \item   The family of discrete conjugate nets
  $\{\xx^\eps\!:B_0^\eps(r)\rightarrow\rz^N\}_{0<\eps<\eps_1}$
  uniquely determined by requiring
  \begin{equation}\label{dcn:goursat}
  x^\eps\rest_{\PP_i^\eps}=X_i\rest_{\PP_i^\eps}\,,\quad
  c_{ij}^\eps\rest_{\PP_{ij}^\eps}=C_{ij}\rest_{\PP_{ij}^\eps}\,,
  \end{equation}
   where $c_{ij}^\eps$ are the rotation coefficients of the net $\xx^\eps$,
   $O(\eps)$-converges in $C^\infty(B_0(r))$ to $\xx$.
  \end{enumerate}
\end{thm}
\noindent
{\bf Proof.} It is easy to reformulate the data (\ref{cn:goursat}),
(\ref{dcn:goursat}) into valid initial data for the Goursat problem
on $\cal X$ for the systems (\ref{eq:ccn1.})--(\ref{eq:ccn3.}),
resp. (\ref{eq:dcn1})--(\ref{eq:dcn3}). Namely, at $\eps=0$ one translates
(\ref{cn:goursat}) into the equivalent data
\begin{equation}\label{cn:goursat alt}
  x(0)=X\,,\quad
  w_i\rest_{\PP_i}=\partial_iX_i\,,\quad c_{ij}\rest_{\PP_{ij}}=C_{ij}\,.
\end{equation}
while for $\eps>0$ one translates (\ref{dcn:goursat}) into the equivalent data
\begin{equation}\label{dcn:goursat alt}
  x^\eps(0)=X\,,\quad
  w_i^\eps\rest_{\PP_i^\eps}=\delta_iX_i\rest_{\PP_i^\eps}\,,\quad
  c_{ij}^\eps\rest_{\PP_{ij}^\eps}=C_{ij}\rest_{\PP_{ij}^\eps}\,.
\end{equation}
Now the theorem follows directly from Lemma \ref{lemma conj}
and Theorem \ref{thm:approx}. \BBox

\begin{thm}\label{thm:jonapprox}
{\bf (approximation of a Jonas pair).}
  Let, in addition to the data listed in Theorem \ref{thm:cnapprox},
  there be given:
  \begin{itemize}
  \item $m$ smooth curves $\XX^+_i:\PP_i\rightarrow\rz^N$,
  $i=1,\ldots,m$, such that all $\XX_i^+$ intersect at a common point
  $\XX^+=\XX^+_1(0)=\ldots=\XX^+_m(0)$,
  and such that for any $1\le i\le m$ and for any point $\xi\in\PP_i$
  the three vectors $\partial_iX_i$, $\partial_iX_i^+$ and
  $\delta_MX_i=X_i^+-X_i$ are coplanar.
  \end{itemize}
  Assume that $\XX^+-\XX$ is not parallel to any of the $m$ vectors
  $\partial_i\XX^+_i(0)$. Define the functions $C_{Mi},C_{iM}:\PP_i\mapsto\rz$
  by the formula
  \begin{equation}\label{jcn:cMi}
  \partial_iX_i^+-\partial_iX_i=C_{Mi}\partial_iX_i+C_{iM}(X_i^+-X_i).
  \end{equation}
  Then, for some $r>0$:
  \begin{enumerate}
  \item
    In addition to the conjugate net $\xx:B_0(r)\rightarrow\rz^N$ defined
    in Theorem \ref{thm:cnapprox}, there is its unique Jonas
    transformation $\xx^+:B_0(r)\rightarrow\rz^N$ such that
    \begin{equation}\label{jcn:goursat}
     x^+\rest_{\PP_i}=X_i^+\quad (1\le i\le m);
    \end{equation}
  \item
    The family $\{\xx^\eps:B^\eps(r)\rightarrow\rz^N\}_{0<\eps<\eps_1}$
    of discrete $M$--dimensional conjugate nets ($M=m+1$), uniquely
    determined by the requirements
    \begin{equation}\label{jdcn:goursat}
    \xx^{\eps}\rest_{\PP_i^\eps}=\XX_i\rest_{\PP_i^\eps}\,,\quad
    c^\eps_{ij}\rest_{\PP_{ij}^\eps}=C_{ij}\rest_{\PP_{ij}^\eps},
    \quad (1\le i\neq j\le m)
    \end{equation}
    and
    \begin{equation}\label{jdcn:goursat+}
    \xx^\eps({\mathbf 0},1)=\XX^+,\quad
    c^\eps_{Mi}\rest_{\PP_i^\eps}=C_{Mi}\rest_{\PP_i^\eps},\quad
    c^\eps_{iM}\rest_{\PP_i^\eps}=C_{iM}\rest_{\PP_i^\eps}
    \quad (1\le i\le m)
    \end{equation}
  $O(\eps)$-converges in $C^\infty(B(r))$ to the Jonas pair $\xx,\xx^+$,
  in the sense that $x^\eps(\cdot,0)\to x$, $x^\eps(\cdot,1)\to x^+$
  \end{enumerate}
\end{thm}
{\bf Proof.} Consider the Goursat problem for $M$--dimensional
discrete conjugate nets with the initial data
\begin{equation}\label{jdcn:goursat alt}
    \xx^{\eps}(0)=\XX\,,\quad
    w^\eps_i\rest_{\PP_i^\eps}=\delta_i\XX_i\rest_{\PP_i^\eps},\quad
    w^\eps_M(0)=\XX^+-\XX, \quad(i=1,\ldots,m)
    \end{equation}
    and
    \begin{equation}\label{jdcn:goursat+ alt}
    c^\eps_{ij}\rest_{\PP_{ij}^\eps}=C_{ij}\rest_{\PP_{ij}^\eps}, \quad
    c^\eps_{Mi}\rest_{\PP_i^\eps}=C_{Mi}\rest_{\PP_i^\eps},\quad
    c^\eps_{iM}\rest_{\PP_i^\eps}=C_{iM}\rest_{\PP_i^\eps}
    \quad (1\le i\neq j\le m).
    \end{equation}
The statement of the theorem follows by applying Theorem
\ref{thm:approx} to this situation. This, in turn, is possible due
to Lemma \ref{lemma jonas} and the following observation. The
condition of the theorem yields that $\XX_i^+(\xi)-\XX_i(\xi)$ is
not parallel to any of $\partial_i\XX^+_i(\xi)$ not only at
$\xi=0$, but also in some neighborhoods of zero on the
corresponding axes $\PP_i$. Now one deduces from the definition
(\ref{jcn:cMi}) of the rotation coefficients $C_{iM}$, $C_{Mi}$
that in the same neighborhoods $C_{Mi}\neq -1$. Hence, the data of
our Goursat problem belong to the set $\DD$ of Lemma \ref{lemma
jonas}. \BBox
\bigskip

Letting not one but two or three directions of the conjugate net
remain discrete in the continuous limit (so that $m'=2$ or $m'=3$), one
arrives at the permutability properties of the Jonas
transformations.
\begin{thm}
{\bf (permutability of Jonas transformations).}
\begin{enumerate}
\item Given an $m$--dimensional
conjugate net $x(\cdot,0,0):B_0(r)\mapsto\rz^N$ and its two Jonas
transformations
$x(\cdot,1,0):B_0(r)\mapsto\rz^N$ and
$x(\cdot,0,1):B_0(r)\mapsto\rz^N$, there exists a
two--parameter family of conjugate nets
$x(\cdot,1,1):B_0(r)\mapsto\rz^N$ that are Jonas transformations
of both $x(\cdot,1,0)$ and $x(\cdot,0,1)$. Corresponding points of
the four conjugate nets are coplanar.
\item
Given three Jonas transformations
\[
x(\cdot,1,0,0),x(\cdot,0,1,0),x(\cdot,0,0,1):B_0(r)\mapsto\rz^N
\]
of a given $m$--dimensional conjugate net
$x(\cdot,0,0,0):B_0(r)\mapsto\rz^N$, as well as three further
conjugate nets
\[
x(\cdot,1,1,0),x(\cdot,0,1,1),x(\cdot,1,0,1):B_0(r)\mapsto\rz^N
\]
such that $x(\cdot,1,1,0)$ is a Jonas transformation of both
$x(\cdot,1,0,0)$ and $x(\cdot,0,1,0)$ etc., there exists
generically a unique conjugate net
\[
x(\cdot,1,1,1):B_0(r)\mapsto\rz^N
\]
which is a Jonas
transformation of all three $x(\cdot,1,1,0),x(\cdot,0,1,1)$ and
$x(\cdot,1,0,1)$.
\end{enumerate}
\end{thm}
These results are proved in the exactly same manner as Theorem
\ref{thm:jonapprox}, if one takes into account the
multi-dimensional consistency of discrete conjugate nets
(Proposition \ref{prop:dcn consistency}).


\section{Orthogonal systems}
\label{sect otho systems}

We give the definitions of orthogonal systems and their
transformations, formulate corresponding Goursat problems, and
prove convergence of discrete orthogonal systems to continuous
ones.

\subsection{Basic definitions}

\begin{dfn} \label{dfn:os}
  A conjugate net $\xx:B_0(r)\mapsto\rz^N$
  is called an $m$-dimensional {\em orthogonal system} in $\rz^N$
  (an {\em orthogonal coordinate system} if $m=N$), if
  \begin{equation}
    \label{eq:corthogonality}
    \partial_i\xx\cdot\partial_j\xx=0.
  \end{equation}
\end{dfn}
\begin{dfn} \label{dfn:rp}
  A pair of $m$--dimensional orthogonal systems $\xx,\xx^+:B_0(r)\mapsto\rz^N$
  is called a {\em Ribaucour pair}, if for any $1\le i\le m$ the
  corresponding coordinate lines of $x$ and $x^+$ envelope a one--parameter
  family of circles, i.e. if it is a Jonas pair of conjugate nets and
  \begin{equation}
    \label{eq:Riba}
    \Big(\frac{\partial_i\xx^+}{|\partial_i\xx^+|}+
    \frac{\partial_i\xx}{|\partial_i\xx|}\Big)
    \cdot(\xx^+-\xx)=0.
  \end{equation}
\end{dfn}
{\bf Remarks.}
\begin{itemize}
\item $\xx^+$ is also called a
Ribaucour transformation of $\xx$; Ribaucour transformations can
be iterated, and our results generalize to any finite sequence of
such transformations (cf. the remarks after Definition
\ref{dfn:jp}).
\item Two-dimensional $(m=2)$ orthogonal systems are
called {\it C-surfaces}. C-surfaces in $\rz^3$  are characterized
as surfaces parametrized along curvature lines.
\item There are different definitions of Ribaucour transformations
in the literature (cf. \cite{udo}). We have chosen the one that is best suited
for our present purposes.
\end{itemize}

\begin{dfn} \label{dfn:dos}
  A map $\xx\!:\!B^\eps(r)\mapsto\rz^N$ is called an $M$-dimensional
  {\em discrete orthogonal system} in $\rz^N$, if the four points
  $\xx$, $\tau_i\xx$, $\tau_j\xx$ and $\tau_i\tau_j\xx$ are
  concircular for all $\xi\in B^\eps(r)$ and all $1\le j\neq j\le M$.

  If $m'=1$, $M=m+1$, then the $M$-dimensional discrete orthogonal
  system, considered as a pair of functions
  $\xx,\xx^+:B^\eps_0(r)\mapsto\rz^N$,
  is called a {\em Ribaucour pair} of $m$-dimensional discrete
  orthogonal systems.
\end{dfn}
{\bf Remark.}
\begin{itemize}
\item Also in the discrete case we shall use the term {\it C-surfaces}
for two-dimensional $(m=2)$ orthogonal systems.
\end{itemize}
So, orthogonal systems form a subclass of conjugate nets subject
to a certain additional condition (orthogonality of coordinate
lines, resp. circularity of elementary quadrilaterals).
\medskip

The classical description of continuous orthogonal systems is this
(in the following equations it is assumed that $i\neq j \neq k \neq i)$:
\begin{eqnarray}
    \partial_i\xx  & = & h_i v_i\,,     \label{eq:DX}\\
    \partial_i v_j & = & \beta_{ji} v_i\,,     \label{eq:DV}\\
    \partial_i h_j & = & h_i\beta_{ij}\,,    \label{eq:DH}\\
    \partial_i\beta_{kj} & = & \beta_{ki}\beta_{ij}\,,    \label{eq:DB}\\
    \partial_i\beta_{ij}+\partial_j\beta_{ji}
    & = & -\partial_i v_i\cdot\partial_j v_j\,.    \label{eq:DL}
\end{eqnarray}
Eq. (\ref{eq:DX}) defines a system $\{v_i\}_{1\le i\leq m}$ of $m$
orthonormal vectors at each point $\xi\in B_0(r)$; the quantities
$h_i=|\partial_i\xx|$ are the respective metric coefficients.
Conjugacy of the net $\xx$ and normalization of $v_i$ imply that
eq. (\ref{eq:DV}) holds with some real--valued functions $\beta_{ij}$.
Eq. (\ref{eq:DH}) expresses the consistency condition
$\partial_i(\partial_j\xx)=\partial_j(\partial_i\xx)$ for $i\neq j$.
Analogously, eq. (\ref{eq:DB}) expresses the consistency condition
$\partial_i(\partial_j v_k)=\partial_j(\partial_i v_k)$ for
$i\neq j\neq k\neq i$. The $m$ equations (\ref{eq:DH}) and $m\choose 3$
equations (\ref{eq:DB}) are called the {\it Darboux system};
they constitute a description of conjugate nets alternative to the system
(\ref{eq:ccn1.})--(\ref{eq:ccn3.}) we used in Sect. \ref{sect conj nets}.
The orthogonality constraint is expressed by $m\choose 2$ additional
equations (\ref{eq:DL}), derived from the identity
$\partial_i\partial_j\langle v_i,v_j \rangle=0$.
In the case $m=N$ the scalar product on the right-hand side
of (\ref{eq:DL}) can be expressed in terms of rotation coefficients
$\beta_{ij}$ only:
\begin{equation}\label{eq:DL2}
  \partial_i\beta_{ij}+\partial_j\beta_{ji}+
  \sum_{k\neq i,j}\beta_{ki}\beta_{kj}=0.
\end{equation}
The Darboux system (\ref{eq:DH}), (\ref{eq:DB}) together with the
eqs. (\ref{eq:DL2}) forms the {\it Lam\'{e} system}.

The classical approach as presented is based on the Euclidean geometry.
However, the invariance group of orthogonal systems is the M\"obius group,
which acts on the compactification $\rz^N\cup\{\infty\}\approx\sx^N$,
rather than on $\rz^N$. This is a motivation to consider orthogonal systems
in the sphere $\sx^N$ and to study their M\"obius-invariant description.
In particular, this enables one to give a frame description of
orthogonal systems which can be generalized to the discrete context
in a straightforward manner. This turns out to be the key to derivation
of the discrete analogue of the Lam\'{e} system.

\subsection{M\"obius geometry}
We give a brief presentation of M\"obius geometry, tailored for our
current needs. A more profound introduction may be found in \cite{udo}.

The $N$--dimensional M\"obius geometry is associated with the
unit sphere $\sx^N\subset\rz^{N+1}$.
Fix two antipodal points on $\sx^N$, $p_0=(0,\ldots,0,1)$
and $p_\infty=-p_0$. The standard stereographic projection $\sigma$
from the point $p_\infty$ is a conformal bijection from
$\sx^N_*=\sx^N\setminus\{p_\infty\}$ to $\rz^N$ that maps spheres
(of any dimension) in $\sx^N$ to spheres or affine subspaces in $\rz^N$.
Its inverse map
\begin{equation*}
  \sigma^{-1}(x)=\frac{2}{1+|x|^2}\,x+\frac{1-|x|^2}{1+|x|^2}\,p_0
\end{equation*}
lifts an orthogonal system in $\rz^N$, continuous or discrete,
to an orthogonal system in $\sx^N$.

The group ${\cal M}(N)$ of $N$-dimensional M\"obius transformations
consists of those bijective maps $\mu:\sx^N\rightarrow\sx^N$ that map
any sphere in $\sx^N$ to a sphere of the same dimension; it then follows
that $\mu$ is also conformal.
\begin{fct}
  The notion of orthogonal system on $\sx^N$ is invariant under M\"obius
  transformations.
\end{fct}
M\"obius transformations can be embedded into a matrix group, namely
the group of pseudo-orthogonal transformations on the
Minkowski space $\rz^{N+1,1}$, i.e. the $(N+2)$-dimensional space spanned
by $\{\ee_1,\ldots,\ee_{N+2}\}$ and equipped with the Lorentz scalar product
\begin{equation*}
  \Big\langle
  \sum_{i=1}^{N+2}u_i\ee_i,\sum_{j=1}^{N+2}v_j\ee_j\Big\rangle =
  \sum_{k=1}^{N+1}u_k v_k - u_{N+2}v_{N+2}.
\end{equation*}
To do this, a model is used where points of $\sx^N$ are identified with
lines on a {\it light cone}
\begin{equation*}
  \cn=\left\{u\in\rz^{N+1,1}\mid\langle u,u \rangle=0\right\}
  \subset \rz^{N+1,1}.
\end{equation*}
This identification is achieved via the projection map
\begin{eqnarray*}
  \pi: \rz^{N+1,1}\setminus(\rz^{N+1}\times\{0\})& \mapsto & \rz^{N+1},\\
  (u_1,\ldots,u_{N+1},u_{N+2}) & \mapsto &
  (u_1/u_{N+2}\,,\ldots,u_{N+1}/u_{N+2}),
\end{eqnarray*}
which sends lines on $\cn$ to points of $\sx^N$. In particular, the
lines through
\begin{equation*}
  \ee_0=\frac{1}{2}(\ee_{N+2}+\ee_{N+1})\quad {\rm and} \quad
  \ee_\infty=\frac{1}{2}(\ee_{N+2}-\ee_{N+1})
\end{equation*}
are mapped to the points $p_0$ and $p_\infty$, respectively.

The group $O^+(N+1,1)$ of genuine Lorentz transformations  consists of
pseudo-orthogonal linear maps $L$ which preserve ``the direction of time'':
\begin{equation}\label{eq:Lorentz1}
  \langle L(u),L(v)\rangle=\langle u,v \rangle\quad
  \forall u,v\in\rz^{N+1,1},\qquad
  \langle L(\ee_{N+2}),\ee_{N+2}\rangle < 0.
\end{equation}
For $L\in O^+(N+1,1)$, the projection $\pi$ induces a M\"obius transformation
of $\sx^N$ via $\mu(L)=\pi\circ L\circ\pi^{-1}$. Indeed, $L$ is linear
and preserves the light cone because of (\ref{eq:Lorentz1}),
so $\mu(L)$ maps $\sx^N$ to itself.
Further, $L$ maps linear planes to linear planes;
under $\pi$, linear planes in $\rz^{N+1,1}$ project to affine
planes in $\rz^{N+1}$. Since any sphere in $\sx^N$ is
uniquely represented as the intersection
of an affine plane in $\rz^{N+1}$ with $\sx^N\subset\rz^{N+1}$,
we conclude that spheres of any dimension are mapped by $\mu(L)$ to
spheres of the same dimension, which is the defining property of M\"obius
transformations. It can be shown that the correspondence $L\mapsto\mu(L)$
between $L\in O^+(N+1,1)$ and ${\cal M}(N)$ is indeed one-to-one.

Recall that we are interested in orthogonal systems in $\rz^N$;
their lifts to $\sx^N$ actually do not contain $p_\infty$. Therefore, we
focus on those M\"obius transformations which preserve $p_\infty$,
thus corresponding to Euclidean motions and homotheties in
$\rz^N=\sigma(\sx^N_*)$. For these, we can restrict our
attention to the section
\begin{equation*}
  {\cal K}=\Big\{u\in\rz^{N+1,1}:\langle u,u\rangle=0,\,
  \langle u,\ee_\infty\rangle=-1/2\Big\}
  =\cn\cap(\ee_0+\ee_\infty^\dashv).
\end{equation*}
(We use here the notation $u^\dashv$ for the hyperplane orthogonal to $u$.)
The canonical lift
\begin{eqnarray*}
  \lambda:\rz^N & \mapsto & {\cal K},\\
  \xx & \mapsto & \xx+\ee_0+|\xx|^2\ee_\infty,
\end{eqnarray*}
factors the (inverse) stereographic projection, $\sigma^{-1}=\pi\circ\lambda$,
and is an isometry in the following sense: for any four points
$\xx_1,\ldots,\xx_4\in\rz^{N}$ one has:
\begin{equation}
  \label{eq:isometry}
  \langle\lambda(\xx_1)-\lambda(\xx_2),\lambda(\xx_3)-\lambda(\xx_4)\rangle
  =(\xx_1-\xx_2)\cdot(\xx_3-\xx_4).
\end{equation}
We summarize these results in a diagram
\begin{displaymath}
  \begin{array}{ccl}
  {\cal K}&\hookrightarrow&\cn\subset\rz^{N+1,1}\\
  {\scriptstyle \lambda}\uparrow&&\downarrow{\scriptstyle \pi}\\
  \rz^N&\xrightarrow{\sigma^{-1}}&\sx^N\subset\rz^{N+1}\\
  \end{array}
\end{displaymath}
Denote by $O_\infty^+(N+1,1)$ the subgroup of Lorentz transformations
fixing the vector $\ee_\infty$ (and consequently, preserving $\cal K$).
The same arguments as before allow one to identify this subgroup with
the group $E(N)$ of Euclidean motions of $\rz^N$.
\begin{fct}
  M\"obius transformations of $\sx^N$ are in a one-to-one correspondence
  with the genuine Lorentz transformations of $\rz^{N+1,1}$.
  Euclidean motions of $\rz^N$ are in a one-to-one correspondence
  with the Lorentz transformations that fix $\ee_\infty$.
\end{fct}
Our next technical device will be Clifford algebras which are very
convenient to describe Lorentz transformations in $\rz^{N+1,1}$, and
hence M\"obius transformations in $\sx^N$.
Recall that the Clifford algebra ${\cal C}(N+1,1)$ is an algebra over $\rz$
with  generators $\ee_1,\ldots,\ee_{N+2}\in\rz^{N+1,1}$ subject to the
relations
\begin{equation}  \label{eq:cliffrule}
  uv+vu=-2\langle u,v\rangle{\mathbf 1}=-2\langle u,v\rangle\qquad
  \forall u,v\in\rz^{N+1,1}.
\end{equation}
Relation (\ref{eq:cliffrule}) implies that $u^2=-\langle u,u \rangle$,
so any vector $u\in\rz^{N+1,1}\setminus\cn$ has an inverse
$u^{-1}=-u/\langle u,u \rangle$. The multiplicative group generated
by the invertible vectors is called the {\it Clifford group}.
It contains the subgroup $Pin(N+1,1)$ which is a universal cover of the
group of M\"obius transformations. We shall need the {\it genuine Pin group}:
\begin{equation*}
  \cH=Pin^+(N+1,1)=\{u_1\cdots u_n\mid u_i^2=-1\},
\end{equation*}
and its subgroup generated by vectors orthogonal to $\ee_\infty$:
\begin{equation*}
  \cH_\infty=\{u_1\cdots u_n\mid
  u_i^2=-1,\;\langle u_i,\ee_\infty\rangle=0\}\subset \cH.
\end{equation*}
$\cH$ and $\cH_\infty$ are Lie groups with Lie algebras
\begin{eqnarray}
  \gh & \!= \!& spin(N+1,1)\;=\;{\rm span}\Big\{\ee_i\ee_j\,:\,
  i,j\in\{0,1,\ldots,N,\infty\},\;i\neq j\Big\},\qquad\\
  \label{eq:Liealgebraarebivectors}
  \gh_\infty & \!=\! & spin_\infty(N+1,1)\;=\;
  {\rm span}\Big\{\ee_i\ee_j\,:\,
  i,j\in\{1,\ldots,N,\infty\},\;i\neq j\Big\}.\quad
\end{eqnarray}
In fact, $\cH$ and $\cH_\infty$ are universal covers
of the previously defined Lorentz subgroups $O^+(N+1,1)$ and
$O^+_\infty(N+1,1)$, respectively. The covers are double, since the
the lift $\psi(L)$ of a Lorentz transformation $L$ is defined up to a sign.
To show this, consider the coadjoint action of $\psi\in \cH$ on
$v\in\rz^{N+1,1}$:
\begin{equation}
  A_\psi(v)=\psi^{-1}v\psi.
\end{equation}
Obviously, for a vector $u$ with $u^2=-1$ one has:
\begin{equation}
  \label{eq:action}
  A_u(v)=u^{-1}vu=2\langle u,v\rangle u-v.
\end{equation}
Thus, $A_u$ is, up to sign, the reflection in the (Minkowski) hyperplane
$u^\dashv$ orthogonal to $u$; these reflections generate the genuine
Lorentz group. The induced M\"obius transformation on $\sx^N$
(for which the minus sign is irrelevant) is the inversion of $\sx^N$
in the hypersphere $\pi(u^\dashv\cap\cn)\subset\sx^N$;
inversions in hyperspheres generate the M\"obius group.
In particular, if $u$ is orthogonal to $\ee_\infty$, then $A_u$ fixes
$\ee_\infty$, hence leaves $\cal K$ invariant and
induces a Euclidean motion on $\rz^N$, namely the reflection in an
affine hyperplane; such reflections generate the Euclidean group.

\begin{fct}
  M\"obius transformations of $\sx^N$ are
  in a one-to-one correspondence to elements of the group $\cH/\{\pm 1\}$;
  Euclidean transformations of $\rz^N$ are
  in a one-to-one correspondence to elements of $\cH_\infty/\{\pm 1\}$.
\end{fct}

\subsection{Orthogonal systems in M\"obius geometry}

We use elements $\psi\in\cH_\infty$ as frames when
describing orthogonal systems in the M\"obius picture.
The vectors $A_\psi(\ee_i)$ for $i=1,\ldots,N$ form an orthogonal
basis in $T{\cal K}$ at the point $A_\psi(\ee_0)\in{\cal K}$.
For an orthogonal system $\xx^\eps: B^\eps(r)\mapsto{\cal K}$,
continuous $(\eps=0)$ or discrete $(\eps>0)$,
with $(m'>0)$ or without $(m'=0)$ transformations,
consider its lift to the conic section $\cal K$:
\begin{equation}
  \hat{\xx}^\eps=\lambda\circ\xx^\eps: B^\eps(r)\rightarrow{\cal K}.
\end{equation}

\subsubsection{Continuous orthogonal systems}

Define $\hat{v}_i:B_0(r)\rightarrow T{\cal K}$ for $i=1,\ldots,m$ by
\begin{equation} \label{eq:MX}
  \partial_i\hat{\xx}=h_i\hat{v}_i,\qquad h_i=|\partial_i \xx|.
\end{equation}
The vectors $\hat{v}_i$ are pairwise orthogonal, orthogonal to $\ee_\infty$,
and can be written as
\begin{equation}\label{hatv}
  \hat{v}_i=v_i+2\langle x,v_i\rangle\ee_\infty
\end{equation}
with the vector fields $v_i$ from (\ref{eq:DX}).
As an immediate consequence, these vectors satisfy
\begin{equation} \label{eq:MV}
  \partial_i\hat{v}_j=\beta_{ji}\hat{v}_i\,,\quad 1\le i\neq j\le m\,,
\end{equation}
with the same rotation coefficients $\beta_{ji}$ as in (\ref{eq:DV}).

\begin{dfn} \label{dfn:framesforcos}
\begin{itemize}
\item[a)]
Given a point $\hat{\xx}\in{\cal K}$ and vectors
$\hat{v}_i\in T_{\hat{x}}{\cal K}$,
$\;1\le i\le m$, we call an element $\psi\in\cH_\infty$ {\em suited} to
$(\hat{\xx};\hat{v}_1,\ldots,\hat{v}_m)$ if
\begin{eqnarray}
    A_\psi(\ee_0) & = & \hat{\xx},    \label{eq:fpoint}\\
    A_\psi(\ee_k) & = & \hat{v}_k\quad(1\leq k\leq m),
    \label{eq:fdirection}
  \end{eqnarray}
\item[b)]
  A frame $\psi:B_0(r)\mapsto\cH_\infty$ is called {\em adapted}
  to the orthogonal system $\hat{\xx}$, if it is suited to
  $(\hat{\xx};\hat{v}_1,\ldots,\hat{v}_m)$ at any point $\xi\in B_0(r)$,
  and
  \begin{equation}\label{eq:fnormal}
   \partial_i A_\psi(\ee_k)=
   \beta_{ki}\hat{v}_i\quad(1\le k\leq N,\;\;1\leq i\leq m).
  \end{equation}
  with some scalar functions $\beta_{ki}$.
\end{itemize}
\end{dfn}

\begin{prp}\label{lem:framesforcos}
  For a given continuous orthogonal system $\hat{x}:B_0(r)\mapsto\cal K$,
  an adapted frame $\psi$ always exists. It is unique, if $m=N$. If $m<N$,
  then an adapted frame is uniquely determined by its value at one point, say
  $\xi=0$. An adapted frame $\psi$ satisfies the following
  differential equations.
  \begin{itemize}
  \item Static frame equations:
    \begin{equation}
      \partial_i\psi=-\ee_i\psi s_i\quad (1\le i\leq m)
    \end{equation}
    with $s_i=(1/2)\partial_i\hat{v}_i$.
  \item Moving frame equations
    \begin{equation} \label{eq:cosmoving}
      \partial_i\psi=-S_i\ee_i\psi\quad (1\le i\leq m)
    \end{equation}
    with
    \begin{equation}  \label{eq:dataforcos}
      S_i=A_{\ee_i\psi}^{-1}(s_i)=\frac{1}{2}
      \sum_{k\neq i}\beta_{ki}\ee_k- h_i\ee_\infty,
    \end{equation}
    where the functions $\beta_{ij}$ and $h_i$ constitute a solution to the
    Lam\'{e} system (in the following equations the indices $1\le i,j\le m$
    and $1\le k\le N$ are pairwise distinct):
    \begin{eqnarray}
      \partial_i h_j & = & h_i\beta_{ij}\,,      \label{eq:ddH}\\
      \partial_i \beta_{kj} & = & \beta_{ki}\beta_{ij}\,,    \label{eq:ddB}\\
      \partial_i \beta_{ij}+\partial_j\beta_{ji} & = &
      -\sum_{k\neq i,j}\beta_{ki}\beta_{kj}\,.      \label{eq:ddL}
    \end{eqnarray}
  \end{itemize}
  Conversely, if $\beta_{ij}$ and $h_i$ are solutions to the
  equations (\ref{eq:ddH})-(\ref{eq:ddL}), then the moving frame
  equations (\ref{eq:cosmoving}) are compatible, and any solution is a
  frame adapted to an orthogonal system.
\end{prp}
{\bf Proof.}
For an orthogonal coordinate system $(m=N)$ the adapted frame is uniquely
determined at any point by the requirements (\ref{eq:fpoint}) and
(\ref{eq:fdirection}), while eqs. (\ref{eq:fnormal}) follow from
(\ref{eq:fdirection}). If $m<N$, extend the $m$ vectors
$\hat{v}_{i}(0)$, $\;1\le i\leq m$, by $N-m$
vectors $\hat{v}_k(0)$, $\;m<k\leq N$, to a orthonormal basis of
$T_{\hat{\xx}(0)}{\cal K}$. There exist unique vector fields
$\hat{v}_k:B_0(r)\mapsto T{\cal K}$, $\;m<k\le N$ with these prescribed
values at $\xi=0$, normalized $(\hat{v}_k^2=-1)$, orthogonal to each
other, to $\hat{v}_i$, $1\le i\le m$, to $\hat{\xx}$ and to $\ee_\infty$,
and satisfying  differential equations
\begin{equation} \label{eq:kevolve}
  \partial_i\hat{v}_k=\beta_{ki}\hat{v}_i\,,\quad
  1\le i\le m,\quad m<k\le N
\end{equation}
with some scalar functions $\beta_{ki}$. Indeed, orthogonality conditions
yield that necessarily
$\beta_{ki}=-\langle\partial_i\hat{v}_i,\hat{v}_k\rangle$.
Eqs. (\ref{eq:kevolve}) with these expressions for $\beta_{ki}$ form a
well-defined linear system of first order partial differential equations
for $\hat{v}_k$. The compatibility condition for this system,
$\partial_j(\partial_i\hat{v}_k)=\partial_j(\partial_i\hat{v}_k)$
$\;(1\le i\neq j\le m)$, is easily verified using (\ref{eq:MV}).
The unique frame $\psi:B_0(r)\mapsto\cH_\infty$ with
$A_{\psi}(\ee_0)=\hat{\xx}$, $\;A_\psi(\ee_k)=\hat{v}_k$ $\;(1\leq k\leq N)$
is adapted to $\hat{\xx}$.

Next, we prove the statements about the moving frame equation.
The equations
\[
\partial_iA_{\psi}(\ee_0)=h_iA_{\psi}(\ee_i),\qquad
\partial_iA_{\psi}(\ee_k)=\beta_{ki}A_{\psi}(\ee_i)
\]
are equivalent to
\[
[\ee_0,(\partial_i\psi)\psi^{-1}]=h_i\ee_i\,,\qquad
[\ee_k,(\partial_i\psi)\psi^{-1}]=\beta_{ki}\ee_i\,.
\]
Since, by (\ref{eq:Liealgebraarebivectors}), $(\partial_i\psi)^{-1}\psi$
is spanned by bivectors, we come to the representation of this element in the
form $(\partial_i\psi)^{-1}\psi=-S_i\ee_i$ with $S_i$ as in
(\ref{eq:dataforcos}). The compatibility conditions for the moving frame
equations, $\partial_i(\partial_j\psi)=\partial_j(\partial_i\psi)$,
are equivalent to
\begin{equation*}
  \partial_j S_i\ee_j+\partial_i S_j\ee_i+
  S_i(\ee_i S_j\ee_i)+S_j(\ee_j S_i\ee_j)=0,
\end{equation*}
which, in turn, are equivalent to the system (\ref{eq:ddH})-(\ref{eq:ddL}).
Conversely, for a solution $\psi:B_0(r)\mapsto\cH_\infty$ of the moving frame
equations, define $\hat{\xx}=A_\psi(\ee_0)$ and $\hat{v}_k=A_\psi(\ee_k)$.
Then $\partial_i\hat{\xx}=h_i\hat{v}_i$, and orthogonality of $\hat{v}_i$
yields that $\hat{\xx}$ is an orthogonal system.

We conclude by showing the static frame equation.
If $\psi$ is adapted, then (\ref{eq:dataforcos}) in combination with the
orthogonality relation of the $v_i$ yields:
\begin{eqnarray*}
  s_i & = & A_{\ee_i\psi}(S_i) \; = \;
  \frac{1}{2}\sum_{k\neq i}\beta_{ki}A_{\ee_i\psi}(\ee_k)
  -h_i A_{\ee_i\psi}(\ee_\infty)
  \\
  & = & -\frac{1}{2}\sum_{k\neq i}\beta_{ki}A_\psi(\ee_k)
  +h_i A_\psi(\ee_\infty)\;=\;
  -\frac{1}{2}\sum_{k\neq i}\beta_{ki}\hat{v}_k+h_i\ee_\infty
  \;=\;\frac{1}{2}\,\partial_i\hat{v}_i\,.
\end{eqnarray*}
The last equality follows from
$\beta_{ki}=-\langle\partial_i\hat{v}_i,\hat{v}_k\rangle$ and the fact that
the $\ee_\infty$--component of $\partial_i\hat{v}_i$ is equal to $2h_i$
(the latter follows from (\ref{hatv})). \BBox

\subsubsection{Discrete orthogonal systems}

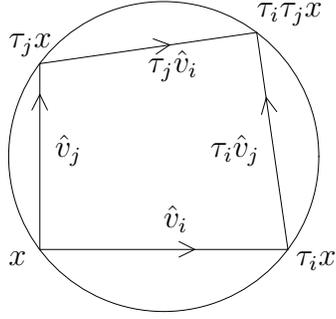
\begin{figure}[htbp]
\begin{center}
\setlength{\unitlength}{0.05em}
\begin{picture}(100,240)(0,-120)
\put(0,0){\circle{200}}
\path(-80,-60)(80,-60)(60,80)(-80,60)(-80,-60)
\path(10,-55)(20,-60)(10,-65)
\path(-85,30)(-80,40)(-75,30)
\path(-7,76)(4,72)(-5,66)
\path(63,27)(66,38)(73,29)
\put(0,-46){$\hat{v}_{i}$}
\put(30,0){$\tau_i\hat{v}_{j}$}
\put(-70,0){$\hat{v}_{j}$}
\put(-10,54){$\tau_j\hat{v}_{i}$}
\put(-100,-70){$\xx$}
\put(85,-70){$\tau_i\xx$}
\put(60,90){$\tau_i\tau_j\xx$}
\put(-100,70){$\tau_j\xx$}
\end{picture}
\caption{An elementary quadrilateral of a discrete orthogonal system}
\label{fig:elem circle}
\end{center}
\end{figure}

Now consider a discrete orthogonal system $\hat{\xx}:B^\eps(r)\mapsto{\cal K}$.
Recall that the cases $m'=0$ and $m'>0$ are treated simultaneously,
so that $\eps_i=\eps$ for $1\le i\leq m$, and $\eps_i=1$ for
$m+1\le i\le m+m'=M$. Appropriate discrete analogues of the
metric coefficients $h_i$ and vectors $\hat{v}_i$ are
\begin{equation}\label{dos:sigmas}
  h_i=|\delta_i\hat{\xx}|,\quad \hat{v}_i=h_i^{-1}\delta_i\hat{\xx}\,.
\end{equation}
These are unit vectors, i.e. $\hat{v}_i^2=-1$, representing the
reflections taking $\hat{\xx}$ to $\tau_i\hat{\xx}$ in the M\"obius picture,
cf. Fig. \ref{fig:elem circle}:
\begin{equation}  \label{eq:xinthemirror}
  \tau_i\hat{\xx}=-A_{\hat{v}_i}\hat{\xx}
  =\hat{\xx}+\eps_i h_i\hat{v}_i\,.
\end{equation}
It is important to note that vectors $\hat{v}_i$ are not mutually orthogonal.
Next, we derive equations that are discrete analogues of (\ref{eq:MV}).
Since four points $\hat{\xx}$, $\tau_i\hat{\xx}$,
$\tau_j\hat{\xx}$ and $\tau_i\tau_j\hat{\xx}$ are coplanar,
the vectors $\tau_j\delta_i\hat{\xx}$ and
$\tau_i\delta_j\hat{\xx}$ lie in the span of
$\delta_i\hat{\xx}$ and $\delta_j\hat{\xx}$.
The lift $\lambda$ is affine, therefore the same holds for $\hat{v}_i$'s:
\begin{equation} \label{eq:sigmasarelinearcombinations}
 \tau_i\hat{v}_j=(n_{ji})^{-1}\cdot(\hat{v}_j+\rho_{ji}\hat{v}_i).
\end{equation}
Here $\rho_{ji}$ are discrete rotation coefficients,
with expected limit behaviour $\rho_{ji}\approx\eps_i\eps_j\beta_{ji}$
as $\eps\to 0$, and
\begin{equation}  \label{eq:ndef}
  n_{ji}^2=1+\rho_{ji}^2+2\rho_{ji}\langle\hat{v}_i,\hat{v}_j\rangle
\end{equation}
is the normalizing factor that is expected to converge to 1 as $\eps\to 0$.
Circularity implies that the angle between $\delta_i\hat{\xx}$ and
$\delta_j\hat{\xx}$ and the angle between
$\tau_j\delta_i\hat{\xx}$ and $\tau_i\delta_j\hat{\xx}$
sum up to $\pi$:
\begin{equation}\label{circularity}
  \langle\tau_j\hat{v}_i,\hat{v}_j\rangle+
  \langle\hat{v}_i,\tau_i\hat{v}_j\rangle=0.
\end{equation}
Under the coplanarity condition (\ref{eq:sigmasarelinearcombinations})
the latter condition is equivalent to
\begin{equation}\label{circ}
  \hat{v}_j(\tau_j\hat{v}_i)+\hat{v}_i(\tau_i\hat{v}_j)=0.
\end{equation}
Inserting (\ref{eq:sigmasarelinearcombinations}) and (\ref{eq:ndef})
into (\ref{circularity}), one finds that either
$\langle\hat{v}_i,\hat{v}_j\rangle=1$, i.e., the circle degenerates to
a line (we exclude this from consideration), or
\begin{equation} \label{eq:dC}
  \rho_{ij}+\rho_{ji}+2\langle\hat{v}_i,\hat{v}_j\rangle=0.
\end{equation}
>From (\ref{eq:ndef}) and (\ref{eq:dC}) there follows:
\begin{equation} \label{eq:thens}
  n_{ij}^2=n_{ji}^2=1-\rho_{ij}\rho_{ji}.
\end{equation}
So, the orthogonality constraint in the discrete case is expressed as the
condition (\ref{eq:dC}) on the rotation coefficients; recall that in
the continuous case this was a condition (\ref{eq:DL}), resp. (\ref{eq:DL2}),
on {\it the derivatives} of the rotation coefficients.

The following definition mimics the static frame equations of the
continuous case.
\begin{dfn}\label{dfn:frames for dos}
  A frame $\psi:B^\eps(r)\mapsto\cH_\infty$ is called
  adapted to a discrete orthogonal system $\hat{\xx}$, if
  \begin{eqnarray}
    A_\psi(\ee_0) & = & \hat{\xx},  \label{eq:pinpointthem}\\
    \tau_i\psi & = & -\ee_i\psi\hat{v}_i\quad (1\le i\leq M).
    \label{eq:pinpointthem1}
  \end{eqnarray}
\end{dfn}
\noindent
Existence of such a frame is a consequence of the circularity condition.
Indeed, the consistency $\tau_j\tau_i\psi=\tau_i\tau_j\psi$ is written as
\[
\ee_i\ee_j\psi\hat{v}_j(\tau_j\hat{v}_i)=
\ee_j\ee_i\psi\hat{v}_i(\tau_i\hat{v}_j).
\]
which is equivalent to (\ref{circ}).
An adapted frame is uniquely defined by the choice of $\psi(0)$ (and
thus is not unique, even if $M=N$). Indeed, eq. (\ref{eq:pinpointthem1})
together with the definition (\ref{dos:sigmas}) of $\hat{v}_i$ imply that
if eq. (\ref{eq:pinpointthem}) holds at one point of $B^\eps(r)$, then it
holds everywhere.

\begin{prp}  \label{lem:framesfordos}
  An adapted frame  $\psi$ for a discrete orthogonal system
  $\hat{\xx}$ satisfies the moving frame equations:
  \begin{equation}  \label{eq:dosmoving}
    \tau_i\psi=-\Sigma_i\ee_i\psi\quad(1\le i\leq M)
  \end{equation}
  with
  \begin{equation} \label{eq:datafordos}
  \Sigma_i=A^{-1}_{\ee_i\psi}(\hat{v}_i)=
  N_i\ee_i+\frac{\eps_i}{2}\sum_{k\neq i}\beta_{ki}\ee_k-\eps_i h_i\ee_\infty,
  \end{equation}
  where
  \begin{equation} \label{eq: N for dos}
  N_i^2=1-\frac{\eps_i^2}{4}\sum_{k\neq i}\beta_{ki}^2\,,
  \end{equation}
  and the functions $\beta_{ij}$ and $h_i$ solve the following {\em discrete
  Lam\'e system} (in the equations below $1\le i,j\leq M$, $1\le k\leq N$,
  and $i\neq j\neq k\neq i$):
  \begin{eqnarray}
      \tau_i h_j & = & (n_{ji}\eps_j)^{-1}\cdot
      (\eps_j h_j+\eps_i h_i\rho_{ij}),    \label{eq:dH}\\
      \tau_i\beta_{kj} & = & (n_{ji}\eps_j)^{-1}\cdot
      (\eps_j\beta_{kj}+\eps_i\beta_{ki}\rho_{ij}),     \label{eq:dB}\\
      \tau_i\beta_{ij} & = & (n_{ji}\eps_j)^{-1}
      \cdot(2N_i\rho_{ij}-\eps_j\beta_{ij}),     \label{eq:dL}\\
      \rho_{ij}+\rho_{ji} & = & N_i\beta_{ij}\eps_j+N_j\beta_{ji}\eps_i
      -\frac{\eps_i\eps_j}{2}\sum_{k\neq i,j}\beta_{ki}\beta_{kj}.
       \label{eq:dQ}
    \end{eqnarray}
    where the abbreviation $n_{ij}=(1-\rho_{ij}\rho_{ji})^{1/2}$
    is used, and $\rho_{ij}$ are suitable real--valued functions.
Conversely, given a solution $\beta_{ij}$ and $h_i$ of the equations
(\ref{eq:dH})-(\ref{eq:dQ}), the system of moving frame equations
(\ref{eq:dosmoving}) is consistent, and its solution $\psi$ is an adapted
frame of a discrete orthogonal system.
\end{prp}
{\bf Proof.} From the definition $\Sigma_i=A^{-1}_{\ee_i\psi}(\hat{v}_i)$
it follows that $\Sigma_i$ is a vector without $\ee_0$--component,
i.e. admits a decomposition of the form (\ref{eq:datafordos}). Again, from
the definition of $\Sigma_i$ and the moving frame equation
(\ref{eq:pinpointthem1}) we find:
\[
A_{\ee_i}(\Sigma_j)=\ee_i\ee_j\psi\hat{v}_j\psi^{-1}\ee_j\ee_i\,,\qquad
\tau_i\Sigma_j=
-\ee_j\ee_i\psi\hat{v}_i(\tau_i\hat{v}_j)\hat{v}_i\psi^{-1}\ee_i\ee_j\,.
\]
Upon using the identity $\hat{v}_i(\hat{v}_j+\rho_{ji}\hat{v}_i)\hat{v}_i=
\hat{v}_j+\rho_{ij}\hat{v}_i\,$, which follows easily from
(\ref{eq:dC}), we can now represent eq. (\ref{eq:sigmasarelinearcombinations})
in the following equivalent form:
\begin{equation} \label{eq:master}
  n_{ji}(\tau_i\Sigma_j)=-A_{\ee_i}(\Sigma_j)-\rho_{ij}A_{\ee_j}(\Sigma_i).
\end{equation}
This equation is equivalent to eqs. (\ref{eq:dH})--(\ref{eq:dL}).
Next, the consistency $\tau_j(\tau_i\psi)=\tau_i(\tau_j\psi)$ of the
moving frame equations (\ref{eq:dosmoving}) is equivalent to
\begin{equation} \label{eq:master1}
  (\tau_i\Sigma_j)A_{\ee_j}(\Sigma_i)+(\tau_j\Sigma_i)A_{\ee_i}(\Sigma_j)=0.
\end{equation}
Inserting the expressions for $\tau_i\Sigma_j$ and $\tau_j\Sigma_i$ from
(\ref{eq:master}) into the left side give:
\begin{eqnarray*}
  A_{\ee_i}(\Sigma_j)A_{\ee_j}(\Sigma_i)+
  A_{\ee_j}(\Sigma_i)A_{\ee_i}(\Sigma_j)+
  \rho_{ji}(A_{\ee_i}(\Sigma_j))^2+
  \rho_{ij}(A_{\ee_j}(\Sigma_i))^2&&\\
  =\;N_i\beta_{ji}\eps_i+N_j\beta_{ij}\eps_j-(\eps_i\eps_j/2)
  \sum_{k\neq i,j}\beta_{ki}\beta_{kj}
  -(\rho_{ij}+\rho_{ji})&=&0,
\end{eqnarray*}
the last equality being nothing but eq. (\ref{eq:dQ}). Conversely,
eqs. (\ref{eq:dH})--(\ref{eq:dQ}) imply (\ref{eq:master}), as well as
the compatibility of the frame equations (\ref{eq:dosmoving}).
So for an arbitrary initial value $\psi(0)\in\cH_\infty$, the frame
equations for $\psi:B^\eps(r)\mapsto\cH_\infty$ can be solved uniquely.
Set $\hat{\xx}^\eps=A_\psi(\ee_0)$, $\hat{v}_i=A_{\ee_i\psi}(\Sigma_i)$.
>From (\ref{eq:master}), (\ref{eq:master1}) for the quantities $\Sigma_i$
there follow (\ref{eq:sigmasarelinearcombinations}), (\ref{circularity})
for the quantities $\hat{v}_i$. Using the fact that $\Sigma_i$ has no
$\ee_0$--component, one shows easily that there holds also
(\ref{eq:xinthemirror}). Therefore, $\hat{x}$ is a discrete orthogonal
system. \BBox

\subsubsection{Ribaucour transformation of a continuous orthogonal \\ system}
To describe a Ribaucour pair $\hat{\xx},\hat{\xx}^+:B_0(r)\mapsto{\cal K}$
of continuous orthogonal systems, one combines the descriptions of
two previous subsections. Recall that in the present context $M=m+1$.

We denote by $h_i^+$, $\hat{v}_i^+$, $\beta_{ij}^+$ the corresponding objects
for $\hat{x}^+$, defined as in (\ref{eq:MX}), (\ref{eq:MV}).
Denote by $\hat{v}_M$ unit vectors such that
\begin{equation}\label{eq: ribdef}
\hat{\xx}^+=-A_{\hat{v}_M}(\hat{\xx}),
\end{equation}
cf. eq. (\ref{eq:xinthemirror}).
The defining property (\ref{eq:Riba}) of Ribaucour transformations and
the normalization of $\hat{v}_M$ imply:
\begin{eqnarray}
  \hat{v}^+_i & = & -A_{\hat{v}_M}(\hat{v}_i),
  \label{eq: ribdef2}\\
  \partial_i\hat{v}_M & = & \frac{\alpha_i}{2}\,(\hat{v}_i^++\hat{v}_i)\;=\;
  \alpha_i(\hat{v}_i-\langle\hat{v}_M,\hat{v}_i\rangle\hat{v}_M),
  \label{eq: ribdef3}
\end{eqnarray}
with $m$ auxiliary functions $\alpha_i:B_0(r)\mapsto\rz$.

\begin{dfn}
  A pair of frames $\psi,\psi^+:B_0(r)\mapsto\cH_\infty$ is called adapted
  to the Ribaucour pair $\hat{\xx},\hat{\xx}^+$, if
  $\psi$ is adapted to $\hat{\xx}$ and
  \begin{equation}
  \psi^+=-\ee_M\psi\hat{v}_M.
  \end{equation}
\end{dfn}

\begin{prp}  \label{lem:framesforrib}
  Let $\psi,\psi^+$ be a pair of frames adapted to the Ribaucour pair
  $\hat{\xx},\hat{\xx}^+$. Then the frame $\psi^+$ is adapted to $\hat{\xx}^+$,
  and the following moving frame equations hold for $\psi^+$:
    \begin{eqnarray}
     \partial_i\psi^+ & = & -S_i^+\ee_i\psi^+\quad(1\le i\leq m),
       \label{eq:ribmoving}\\
     \psi^+ & = & -\Sigma_M\ee_M\psi, \label{eq:ribmoving1}
    \end{eqnarray}
    where
    \begin{eqnarray}
      S_i^+ & = & \frac{1}{2}\sum_{k\neq i}\beta^+_{ki}\ee_k-
      h^+_i\ee_\infty\quad(1\le i\leq m),      \label{eq:dataforrib1}\\
      \Sigma_M & = & A^{-1}_{\ee_M\psi}(\hat{v}_M)\;=\;
      N_M\ee_M+\frac{1}{2}\sum_{k\neq M}\beta_{kM}\ee_k-
      h_M\ee_\infty,       \label{eq:dataforrib2}\\
      N_M^2 & = & 1-\frac{1}{4}\sum_{k\neq M}\beta_{kM}^2,
    \end{eqnarray}
    the functions $h_i,\beta_{ij}$ and $h_i^+,\beta_{ij}^+$ solve eqs.
    (\ref{eq:ddH})-(\ref{eq:ddL}),
    and the following system is satisfied (for
    $1\le i\neq j\leq m$, $\;1\le k\leq N$, and $\;i\neq k\neq M\neq i$):
      \begin{eqnarray}
        \label{eq:ddd1}
        h^+_i&=&h_i+h_M\alpha_i\,,\\
    \label{eq:ddd2}
        \beta^+_{ki}&=&\beta_{ki}+\beta_{kM}\alpha_i\,,\\
        \label{eq:ddd3}
        \beta_{Mi}^+&=&-\beta_{Mi}+2N_M \alpha_i\,,\\
        \label{eq:ddd4}
        \partial_i h_M&=&\frac{1}{2}\,(h_i+h_i^+)\beta_{iM}\,,\\
    \label{eq:ddd5}
        \partial_i \beta_{kM}&=&
    \frac{1}{2}\,(\beta_{ki}+\beta_{ki}^+)\beta_{iM}\,,\\
        \label{eq:ddd8}
        \partial_i\beta_{iM}&=&-\frac{1}{2}\sum_{k\neq i,M}
    (\beta_{ki}+\beta^+_{ki})\beta_{kM}+N_M(\beta_{Mi}-\beta_{Mi}^+)\,,\\
        \label{eq:ddd9}
        \partial_i \alpha_j&=&\frac{1}{2}\,\alpha_i(\beta_{ij}+\beta^+_{ij})\,.
      \end{eqnarray}
  Conversely, given solutions $h_i,\beta_{ij}$ and $h_i^+,\beta_{ij}^+$
  of the system above with suitable auxiliary functions $\alpha_i$,
  the moving frame equations (\ref{eq:ribmoving}), (\ref{eq:ribmoving1})
  are compatible, and the solution $\psi,\psi^+$ is a pair of frames
  adapted to a Ribaucour pair of orthogonal systems.
\end{prp}
{\bf Proof.} First we show that $\psi^+$ is an adapted frame for $\xx^+$.
We have for $1\le i\leq m$, $1\le k\le N$:
\begin{eqnarray*}
  A_{\psi^+}(\ee_0) & = & A_{\hat{v}_M}A_{\ee_M\psi}(\ee_0)=
  -A_{\hat{v}_M}(\hat{\xx})=\hat{\xx}^+\,,\\
  A_{\psi^+}(\ee_i) & = & A_{\hat{v}_M}A_{\ee_M\psi}(\ee_i)=
  -A_{\hat{v}_M}(\hat{v}_i)=\hat{v}^+_i\,,\\
  \partial_i A_{\psi^+}(\ee_k) & = &
  \partial_i(A_{\hat{v}_M}A_{\ee_M\psi}(\ee_k))\;=\;
  \pm\partial_i(A_{\hat{v}_M}A_{\psi}(\ee_k))\\
   & = & \pm\partial_i(A_{\psi}(\ee_k)
  -2\langle\hat{v}_M,A_{\psi}(\ee_k)\rangle\hat{v}_M)\\
  & = &\pm(\beta_{ki}-2\alpha_i\langle\hat{v}_M,A_{\psi}(\ee_k)\rangle)
  (\hat{v}_i-2\langle\hat{v}_M,\hat{v}_i\rangle\hat{v}_M)\\
  & = & \pm (\beta_{ki}-2\alpha_i\langle\hat{v}_M,A_{\psi}(\ee_k)\rangle)
  \hat{v}^+_i\,=\,\beta_{ki}^+\hat{v}^+_i,
\end{eqnarray*}
where the minus sign applies iff $k=M$ (and if $k=0$, the latter calculation
goes through almost literally, with replacing $\beta_{ki}$, $\beta_{ki}^+$
by $h_i$, $h^+$, respectively). Next, from
$\Sigma_M=A^{-1}_{\ee_M\psi}(\hat{v}_M)$ there follows:
\begin{eqnarray*}
  \beta_{kM} & = & 2\langle\Sigma_M,\ee_k\rangle=
    2\langle \hat{v}_M,A_{\ee_M\psi}(\ee_k)\rangle
    =-2\langle\hat{v}_M,A_\psi(\ee_k)\rangle,\\
  h_M & = & 2\langle\Sigma_M,\ee_0\rangle=
    2\langle \hat{v}_M,A_{\ee_M\psi}(\ee_0)\rangle
    =-2\langle\hat{v}_M,\xx\rangle,\\
  N_M & = & \langle\Sigma_M,\ee_M\rangle=
    \langle \hat{v}_M,A_{\ee_M\psi}(\ee_M)\rangle
    =\langle\hat{v}_M,A_\psi(\ee_M)\rangle.
\end{eqnarray*}
This proves eqs. (\ref{eq:ddd1})--(\ref{eq:ddd3}). Further, eqs.
(\ref{eq:ddd4})--(\ref{eq:ddd8}) are readily derived by calculating
the $i$-th paartial derivative of the respective scalar product.
Finally, eq. (\ref{eq:ddd9}) comes from the consistency condition
$\partial_i(\partial_j\hat{v}_M)=\partial_j(\partial_i\hat{v}_M)$.

Conversely, given a solution to the equations (\ref{eq:ddH})--(\ref{eq:ddL})
and (\ref{eq:ddd1})--(\ref{eq:ddd9}), the moving frame equations are
consistent, thus defining the frames $\psi$, $\psi^+$.
It follows from Proposition \ref{lem:framesforcos} that both
$\hat{\xx}=A_{\psi}(\ee_0)$ and $\hat{\xx}^+=A_{\psi^+}(\ee_0)$
are orthogonal systems. Furthermore, eq. (\ref{eq:ribmoving1}),
yields the defining relations
(\ref{eq: ribdef})--(\ref{eq: ribdef3}) of the Ribaucour pair,
with $\hat{v}_M=A_{\ee_M\psi}(\Sigma_M)$.
\BBox


\section{Goursat problems and approximation for \\
orthogonal systems}
\label{sect: ortho approx}

It would be tempting to derive the theory of the Lam\'e system
(\ref{eq:DX})-(\ref{eq:DL}) from its discrete counterpart
(\ref{eq:dH})-(\ref{eq:dQ}), treating the latter as a hyperbolic system
of the type (\ref{eq:evolution}). However, it turns out that in dimensions
$M\geq 3$ this is hard to carry out, since one needs to enlarge the set
of dependent variables and equations in a cumbersome manner.
The way around is based on the following fundamental lemma,
which allows one to take care of two--dimensional orthogonal systems (in
the coordinate planes) only, and to use then the results for conjugate
systems.

\begin{lem}  \label{lem:osaresubcn}
  a) If for a conjugate net $\xx$, continuous or discrete,
  its restriction to each plane $\PP_{ij}$, $1\leq i\neq j\leq M$,
  is a C-surface, then $\xx$ is a (continuous or discrete) orthogonal
  system.

  b) If in a Jonas pair $(\xx,\xx^+)$, continuous or discrete,
  the net $\xx$ is an orthogonal system, and the corresponding
  coordinate curves $x\rest_{\PP_i}$ and $x^+\rest_{\PP_i}$ envelope
  one--dimensional families of circles, then $(\xx,\xx^+)$ is a Ribaucour
  pair of orthogonal systems.
\end{lem}
{\bf Proof} in the discrete case is based on the Miguel theorem, cf.
Sect. \ref{sect:intro}, and can be found in \cite{DoSa1}. For the
continuous case it is a by-product of the proof of Theorem
\ref{thm:dosapprox} below. \BBox
\medskip

So, suppose that $M=2$. Then the moving frame equations
(\ref{eq:dosmoving}) and the system (\ref{eq:dH})-(\ref{eq:dQ})
take the form ($1\le i,j\le 2$, $\;1\le k\le N$, $\;i\neq j\ne
k\ne i$):
\begin{eqnarray}
  \delta_i\psi & = &
  \bigg(\frac{N_i-1}{\eps_i}-\frac{1}{2}\,\sum_{k\neq i}\beta_{ki}\ee_k\ee_i
  +h_i\ee_\infty\ee_i\bigg)\psi\,,  \label{eq:dP.}\\
  \delta_i h_j & = & \frac{\rho_{ij}}{\eps_jn}\,h_i
  +\frac{1-n}{\eps_in}\,h_j\,,  \label{eq:dH.}\\
  \delta_i\beta_{ij} & = & \frac{2N_i\rho_{ij}}{\eps_i\eps_jn}
  -\frac{1+n}{\eps_in}\,\beta_{ij},  \label{eq:dL.}\\
  \delta_i\beta_{kj} & = & \frac{1-n}{\eps_in}\,\beta_{kj}
      +\frac{\rho_{ij}}{\eps_jn}\,\beta_{ki},     \label{eq:dB.}\\
  \rho_{12}+\rho_{21} & = &
  \eps_2 N_1\beta_{12}+\eps_1 N_2\beta_{21}-\eps_1\eps_2\Theta.
  \label{eq:tworos}
\end{eqnarray}
Here we use the abbreviations
\begin{equation}
n=n_{12}=n_{21}=\sqrt{1-\rho_{12}\rho_{21}},\quad
\Theta=\frac{1}{2}\sum_{k>2}\beta_{k1}\beta_{k2},\quad
N_i^2=1-\frac{\eps_i^2}{4}\sum_{k\neq i}\beta_{ki}^2\,.
\end{equation}
We show how to re-formulate the above system in the hyperbolic
form and to pose a Cauchy problem for it.

\subsection{Hyperbolic equations and approximation for \\ C-surfaces}

First consider the case $m=2$, $m'=0$, $\eps_1=\eps_2=\eps$,
corresponding to C-surfaces. Set
\begin{eqnarray}
  \rho_{12} & = & \eps N_1\beta_{12}-(\eps^2/2)(\Theta-\gamma),
    \label{eq:resolve2}\\
  \rho_{21} & = & \eps N_2\beta_{21}-(\eps^2/2)(\Theta+\gamma),
    \label{eq:resolve}
\end{eqnarray}
with a suitable function $\gamma:\PP_{12}^\eps\mapsto\rz$. It can be
said that the auxiliary function $\gamma$ splits the constraint, making
the Lam\'{e} system hyperbolic. This splitting plays a crucial role for our
approximation results. In the smooth limit
$2\gamma=\partial_1\beta_{12}-\partial_2\beta_{21}$ (cf. eqs.
(\ref{eq:DLI}), (\ref{eq:DLII}) below).

As a Banach space for the system take
\[
  \BX=\hat{\cH}\{\psi\}\times\rz^{2(N-1)}\{\beta\}\times
  \rz^2\{h_1,h_2\}\times\rz\{\gamma\}\,,
\]
where $\hat{\cH}$ is the space of an (arbitrary) matrix
representation of $Spin(N+1,1)$ (note that if
$\psi\in\cH_\infty\subset\hat{\cH}$ at some point, then eq.
(\ref{eq:dP.}) guarantees that this is the case everywhere),
$\beta$ denotes the collection of $\beta_{kj}$ with $k=1,2$,
 $1\le j\le N$, $k\neq j$. Obviously,
\begin{equation}\label{Nas}
N_1=1+O(\eps^2),\qquad N_2=1+O(\eps^2),\qquad n=1+O(\eps^2),
\end{equation}
where the constants in $O$--symbols are uniform on compact subsets
of $\BX$.

Both directions $i=1,2$ are assumed to be evolution directions for
$\psi$; for $h_i$ and ${\beta}_{ki}$ there is one evolution
direction $j\ne i$, while for the additional function $\gamma$
both directions $i=1,2$ are static.
\begin{gou} {\bf (for discrete C-Surfaces).} \label{gou:dos}
  Given $\Psi^\eps\in\cH_\infty$, functions $H^\eps_i:\PP^\eps_i\mapsto\rz$
  and $B^\eps_{ki}:\PP^\eps_i\mapsto\rz$ for
  $i=1,2$, $\;1\le k\le N$, $\;k\neq i$, and
  $\Gamma^\eps:\PP^\eps_{12}\mapsto\rz$,
  find a solution to the equations (\ref{eq:dP.})-(\ref{eq:dB.})
  with (\ref{eq:resolve2}), (\ref{eq:resolve}), and
  with the Goursat data
\[
  \psi(0)=\Psi^\eps,\quad h_i\rest_{\PP_i^\eps}=H^\eps_i,
  \quad \beta_{ki}\rest_{\PP_i^\eps}=B^\eps_{ki},
  \quad \gamma\rest_{\PP_{12}^\eps}=\Gamma^\eps.
\]
\end{gou}
\begin{lem} \label{lem:dos}
  The hyperbolic system (\ref{eq:dP.})--(\ref{eq:dB.})
  with (\ref{eq:resolve2}), (\ref{eq:resolve}) is consistent.
  Goursat problem \ref{gou:dos} for discrete C-surfaces satisfies
  condition (F). The limiting system for $\eps=0$ consists of equations
  (\ref{eq:cosmoving}) with (\ref{eq:dataforcos}), (\ref{eq:ddH}),
  (\ref{eq:ddB}) for paiwise distinct indices $1\le i,j\le 2$, $\;1\le k\le N$,
  and
  \begin{eqnarray}
      \partial_1 \beta_{12} & = &
      -(1/2)\sum_{k>2}\beta_{k1}\beta_{k2}+\gamma\,, \label{eq:DLI}\\
      \partial_2 \beta_{21} & = &
      -(1/2)\sum_{k>2}\beta_{k1}\beta_{k2}-\gamma\,, \label{eq:DLII}
    \end{eqnarray}
   which describe continuous C-surfaces.
\end{lem}
{\bf Proof.} Consistency is easy to see: the only condition to be
checked is $\delta_1\delta_2\psi=\delta_2\delta_1\psi$, but this
has already been shown in Proposition \ref{lem:framesfordos}. Also other
statements are obvious from (\ref{Nas}). Notice that hyperbolic equations
(\ref{eq:DLI}), (\ref{eq:DLII}) come to replace the non--hyperbolic
orthogonality constraint (\ref{eq:ddL}). \BBox
\medskip

Now we discuss the way to prescribe the data in the Goursat problem
\ref{gou:dos} in order to get an approximation of a given smooth
C--surface. Unlike in the case of general conjugate nets, it is
not possible to prescribe the discrete curves
$\XX^\eps_i=\xx^\eps\rest_{\PP^\eps_i}$ coinciding with their continuous
counterparts $\XX_i=\xx\rest_{\PP_i}$ at the lattice points.
However, it is still possible to achieve that $\XX^\eps_i$ are completely
determined by $\XX_i$.

Given a curve $\hat{X}_i:\PP_i\mapsto{\cal K}$, one has the corresponding
tangential vector field $\partial_i\hat{X}_i$, and therefore the function
$h_i=|\partial_i\hat{X}_i|$ and the field of unit vectors
$\hat{v}_i=h_i^{-1}\partial_i\hat{X}_i$. Take
$\Psi\in\cH_\infty$ suited to $(\hat{\XX}_i(0);\hat{v}_i(0))$.
Then by Proposition \ref{lem:framesforcos}, there is a unique frame
$\psi:\PP_i\rightarrow\cH_\infty$ adapted to the curve
$\hat{\XX}_i$ with $\psi(0)=\Psi$; this frame is defined as the unique
solution of the differential equation $\partial_i\psi=-S_i\ee_i\psi$ with
$S_i$ given by eq. (\ref{eq:dataforcos}). According to the proof of Proposition
\ref{lem:framesforcos}, the latter equation is equivalent to the system of
equations
\[
\partial_i\hat{v}_k=\beta_{ki}\hat{v}_i, \quad
\beta_{ki}=-\langle \hat{v}_k,\partial_i\hat{v}_i\rangle,
\quad k\neq i.
\]
So, in order to determine the rotation coefficients
$\beta_{ki}:\PP_i\mapsto\rz$ for $1\le k\le N$, $\;k\neq i$,
one has to solve the latter system of ordinary differential
equations with the initial data $\hat{v}_k(0)=A_{\Psi}(\ee_k)$.
Thus, we produced the functions $h_i$ and $\beta_{ki}$, or, what is
equivalent, the Clifford elements
\[
S_i=\frac{1}{2}\sum_{k\neq i}\beta_{ki}\ee_k-h_i\ee_{\infty}
\]
for a given curve $\hat{X}_i:\PP_i\mapsto\cal K$. We say that $h_i$,
$\beta_{ki}$ are {\it read off} the curve $\hat{X}_i$.

\begin{dfn}
 The {\em canonical discretization of the curve
 $\hat{\XX}_i:\PP_i\mapsto\cal K$
 with respect to the initial frame} $\Psi\in\cH_\infty$ suited to
 $(\hat{X}_i(0);\hat{v}_i(0))$ is the function
 $\hat{\XX}^\eps_i=A_{\psi^{\eps}}(\ee_0):\PP^\eps_i\mapsto\cal K$,
 where $\psi^\eps:\PP^\eps_i\mapsto\cH_\infty$ is the solution of the
 discrete moving frame equation $\tau_i\psi^\eps=-\Sigma_i\ee_i\psi^\eps$
 with
\[
\Sigma_i=\bigg(1-\frac{\eps^2}{4}\sum_{k\neq i}\beta_{ki}^2\bigg)^{1/2}\ee_i
+\frac{\eps}{2}\sum_{k\neq i}\beta_{ki}\ee_k-\eps h_i\ee_{\infty}
\]
and the initial condition $\psi^{\eps}(0)=\Psi$.
\end{dfn}
In other words, for the canonical discretization the data $h_i$ and
$\beta_{ki}$ are read off the continuous curve. Obviously, $\psi^\eps$
is an adapted frame for the discrete curve $\XX^\eps_i$. As $\eps\to 0$,
the canonical discretization $\XX_i^\eps$ converges to $\XX_i$ with the rate
$O(\eps)$ in $C^\infty$.

\begin{prp} \label{prp:clsdoexist}
{\bf (approximation for C--surfaces).}
  Let $m=M=2$, and let there be given:
  \begin{itemize}
  \item two smooth curves $\XX_i:\PP_i\mapsto\rz^N$ ($i=1,2$),
  intersecting orthogonally at $X=\XX_1(0)=\XX_2(0)$,
  \item a smooth function $\Gamma:\PP_{12}\mapsto\rz$.
  \end{itemize}
  Assume $\Psi\in\cH_\infty$ is suited for
  $(\hat{\XX};\hat{v}_1(0),\hat{v}_2(0))$. Then, for some $r>0$:
  \begin{enumerate}
  \item There exists a unique C--surface $\xx:B_0(r)\mapsto\rz^N$
    that coincides with $\XX_i$ on $\PP_i\cap B_0(r)$ and satisfies
    \begin{equation*}
      \partial_1\beta_{12}-\partial_2\beta_{21}=2\Gamma.
    \end{equation*}
  \item
    Consider the family of discrete C--surfaces
    $\{\xx^\eps:B_0^\eps(r)\rightarrow\rz^N\}_{0<\eps<\eps_1}$
    defined as the solutions to the Goursat problem \ref{gou:dos} with
    the data
     \[
     \Psi^\eps=\Psi,\quad H^\eps_i=H_i\rest_{\PP_i^\eps},
     \quad B^\eps_{ki}=B_{ki}\rest_{\PP_i^\eps},
     \quad \Gamma^\eps=\Gamma\rest_{\PP_{12}^\eps},
     \]
    where the functions $H_i$ and $B_{ki}$ are read off the smooth curves
    $X_i$, so that $\xx^\eps\rest_{\PP^\eps_i}$ are the canonical
    $\eps$-discretizations of $\XX_i$. This family of discrete C--surfaces
    $O(\eps)$-converges in $C^\infty$ to $\xx$.
  \end{enumerate}
\end{prp}
\noindent
{\bf Proof} follows from Lemma \ref{lem:dos} and Theorem \ref{thm:approx}.
\BBox

\begin{figure}[ht]
\begin{center}
\epsfig{file=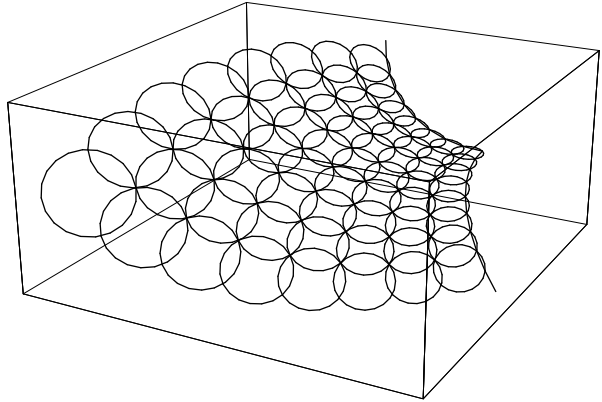,height=150pt}
\end{center}
\caption{Two curves determine a discrete C--surface
via Proposition \ref{prp:clsdoexist}.}
\end{figure}

\subsection{Goursat problem and approximation for an \\ orthogonal system}

\begin{thm} \label{thm:dosapprox}
{\bf (approximation of an orthogonal system).}
  Let $m(m-1)/2$ smooth C--surfaces ${\cal S}_{ij}:\PP_{ij}\mapsto\rz^N$
  be given, labelled by $1\leq i\neq j\leq m$ with
  ${\cal S}_{ij}={\cal S}_{ji}$.
  Assume that for a given $1\le i\le m$ all surfaces ${\cal S}_{ij}$
  intersect along the curvature lines $\XX_i={\cal S}_{ij}\rest_{\PP_i}$.
  Set $h_i=|\partial_i\XX_i|$, assume that $h_i(0)\neq 0$, and set further
  $v_i=h_i^{-1}\partial_iX_i$. All curves $X_i$ intersect at one point
  $\XX=\XX_1(0)=\ldots=\XX_m(0)$ orthogonally.
  Let $\Psi\in\cH_\infty(N)$ be suited to
  $(\hat{\XX},\hat{v}_1(0),\ldots,\hat{v}_m(0))$.
  Construct the discrete C--surfaces ${\cal S}_{ij}^\eps:
  \PP_{ij}^\eps\mapsto\rz^N$ according to Proposition \ref{prp:clsdoexist},
  with the functions
  $\Gamma_{ij}=(\partial_i\beta_{ij}-\partial_j\beta_{ji})/2$ for
  $1\le i<j\le m$. Then for some $r>0$:
  \begin{enumerate}
  \item
    There exists a unique orthogonal system $\xx:B_0(r)\mapsto\rz^N$,
    coinciding with ${\cal S}_{ij}$ on $\PP_{ij}\cap B_0(r)$.
  \item
    There exists a unique family of discrete orthogonal systems
    $\{\xx^\eps\!:\!B_0^\eps(r)\!\mapsto\rz^N\}_{0<\eps<\eps_0}$
    coinciding with ${\cal S}_{ij}^\eps$ on
    $\PP_{ij}^\eps\cap B_0^\eps(r)$.
    The family $x^\eps$ converges to $x$ with the rate
    $O(\eps)$ in $C^\infty$.
  \end{enumerate}
\end{thm}
\noindent
{\bf Proof.}
C--surfaces are special two-dimensional conjugate nets, therefore the
surfaces ${\cal S}_{ij}$ can be supplied with the respective coefficients
$C_{ij},C_{ji}:\PP_{ij}\mapsto\rz$. Now Theorem \ref{thm:cnapprox}
can be applied. It yields the existence of a unique conjugate net
$\xx$ which has $\XX_i$ as image of the $i$-th coordinate axes and
$C_{ij},C_{ji}$ as coefficients on the respective $\PP_{ij}$.

Similarly, discrete C--surfaces are special discrete two-dimensional
conjugate nets, so for each ${\cal S}^\eps_{ij}$, the
coefficients $C^\eps_{ij},C^\eps_{ji}$ are defined,
and they converge in $C^\infty$ to the coefficients
$C_{ij},C_{ji}$ of the respective ${\cal S}_{ij}$.
Note that for any $1\le i\le m$ the discrete surfaces ${\cal S}^\eps_{ij}$
intersect along the discrete curve
$\XX^\eps_i={\cal S}^\eps_{ij}\rest_{\PP_i}$, which is the canonical
discretization of the curve $\XX_i$. The
conclusion 2 of Theorem \ref{thm:cnapprox} is still valid if the Goursat
data for $c_{ij}^\eps$ in (\ref{dcn:goursat}) are taken as $C_{ij}^\eps$
instead of $C_{ij}$, i.e. read off ${\cal S}_{ij}^\eps$ rather than
${\cal S}_{ij}$. Thus, Theorem \ref{thm:cnapprox} delivers a family
of discrete conjugate nets $\{\xx^\eps\}$ that is $O(\eps)$-convergent
to $x$. Since the restrictions of $\xx^\eps$ to the coordinate
planes $\PP^\eps_{ij}$ are discrete C--surfaces, the discrete variant
of Lemma \ref{lem:osaresubcn} implies that the nets $\xx^\eps$ are
actually discrete orthogonal systems. It remains to show that the limiting
net $x$ is with necessity an orthogonal one. But this follows immediately
from the fact of $C^\infty$--convergence and eq. (\ref{circ}): indeed, we
can conclude that for the net $x$ everywhere holds
\[
  \hat{v}_j\hat{v}_i+\hat{v}_i\hat{v}_j=0\quad\Leftrightarrow\quad
  \langle \hat{v}_i,\hat{v}_j\rangle=0,
\]
that is, $\partial_ix\cdot\partial_jx=0$. \BBox
\medskip

This result immediately implies Theorem \ref{thm:cDupin} from Introduction:
to bring its assumptions into the form of Theorem \ref{thm:dosapprox} above,
one needs only to introduce a curvature line parametrization on surfaces
${\cal F}_i$ such that the curves $\XX_i$ are parametrized by arc-length.

\subsection{Goursat problems and approximation for \\
Ribaucour transformations}

\begin{figure}[ht]
\begin{center}
\epsfig{file=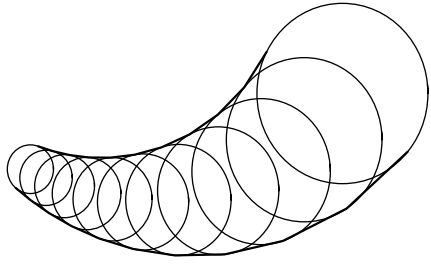,height=100pt}
\end{center}
\caption{A pair of curves enveloping a circle congruence}
\label{fig:envelope}
\end{figure}

Now we consider again a two--dimensional ($M=2$) discrete orthogonal system
described by eqs. (\ref{eq:dP.})--(\ref{eq:dB.}), but perform a different
continuous limit, leading to two curves enveloping a circle congruence
(Fig. \ref{fig:envelope}).
In other words, we set $m=m'=1$, $\eps_1=\eps$, $\eps_2=1$. Recall that in
this case
\[
B^\eps(r)=B_0^\eps(r)\times\{0,1\},\quad {\rm where}\quad
B_0^\eps(r)=\PP_1^\eps\cap [0,r]=(\eps\gz)\cap [0,r],
\]
so that only two lines in $\PP_{12}^\eps$ are considered. Therefore,
a function on $\PP_{12}^\eps$ is conveniently considered as a pair of
functions on $\PP_1^\eps$, or even as just one function on $\PP_1^\eps$,
if one is interested in $B_0^\eps(r)\times\{0\}$ only (this will be the
case for the function $\alpha$ below). We denote the shift $\tau_2$ by
the superscript ``$+$''.

Set
\begin{equation}  \label{eq:resolve1}
  \rho_{21}=\eps\alpha\,,\qquad
  \rho_{12}=N_1\beta_{12}+\eps(N_2\beta_{21}-\Theta-\alpha);
\end{equation}
with a suitable function $\alpha:\PP^\eps_1\mapsto\rz$.

As the Banach space where the system lives we choose, exactly as before,
\[
  \BX=\hat{\cH}\{\psi\}\times\rz^{2(N-1)}\{\beta\}\times
  \rz^2\{h_1,h_2\}\times\rz\{\alpha\}\,.
\]
In the present case we have:
\begin{equation}  \label{eq:Npositive}
  N_1^2=1+O(\eps^2),  \quad
  N_2^2=1-\frac{1}{4}\sum_{k\neq 2}\beta_{k2}^2\,,\quad
  n=1-\frac{\eps}{2}\,\alpha\beta_{12}+O(\eps^2)\,,
\end{equation}
where the constants in $O$--symbols are uniform on compact subsets
of $\BX$. However, the system itself is now well defined not on all of
$\BX$ but rather on the subset where $N_2^2>0$:
\begin{equation}\label{dRib:D}
\DD=\Big\{(\psi,\beta,h_1,h_2,\alpha):\sum_{k\neq 2}
 \beta_{k2}^2<4\Big\}\subset\BX.
\end{equation}
\begin{gou} \label{gou:dos1}
{\bf (for a pair of curves enveloping a circle congruence).}
  Given $\Psi^\eps\in\cH_\infty$, functions $H^\eps_1:\PP^\eps_1\mapsto\rz$
  and $B^\eps_{k1}:\PP^\eps_1\mapsto\rz$ for $2\le k\le N$, and
  $A^\eps:\PP^\eps_1\mapsto\rz$, as well as the real numbers
  $H_2^\eps$ and $B^\eps_{k2}$ for $1\le k\le N$, $\;k\neq 2$,
  find a solution to the equations (\ref{eq:dP.})-(\ref{eq:dB.})
  with (\ref{eq:resolve1}) on $B^\eps(r)$ with the Goursat data
\begin{eqnarray*}
 & \psi(0)=\Psi^\eps,\quad h_1\rest_{\PP_1^\eps}=H^\eps_1,
  \quad \beta_{k1}\rest_{\PP_1^\eps}=B^\eps_{k1},
  \quad \alpha\rest_{\PP_1^\eps}=A^\eps, & \\
 & h_2(0)=H_2^\eps,\quad \beta_{k2}(0)=B_{k2}^\eps. &
\end{eqnarray*}
\end{gou}
\begin{lem}\label{lem:dos1}
  The hyperbolic system (\ref{eq:dP.})--(\ref{eq:dB.})
  with (\ref{eq:resolve1}) is consistent. \linebreak
  Goursat problem \ref{gou:dos1}
  satisfies condition (F) on the subset $\DD$. The limiting system for
  $\eps=0$ consists of equations (\ref{eq:cosmoving}) with
  (\ref{eq:dataforcos}) for $i=1$, (\ref{eq:dosmoving}) with
  (\ref{eq:datafordos}) for $i=2$, and
  \begin{eqnarray}
   \partial_1 h_2 & = & \beta_{12}(h_1+\alpha h_2/2)\,,
   \label{riblemma 1}\\
   \partial_1\beta_{k2} & = & \beta_{12}(\beta_{k1}+\alpha\beta_{k2}/2)\,,
   \quad k>2\,,  \label{riblemma 2}\\
   \partial_1\beta_{12} & = &
   2(N_2\beta_{21}-\Theta-\alpha)+\alpha\beta_{12}^2/2\,,\label{riblemma 3}\\
   h_1^+-h_1 & = & \alpha h_2\,,   \label{riblemma 4}  \\
   \beta_{k1}^+-\beta_{k1} & = & \alpha\beta_{k2}\,,\quad k>2\,,
   \label{riblemma 5}\\
   \beta_{21}^+-\beta_{21} & = & 2(N_2\alpha-\beta_{21})\,,
   \label{riblemma 6}
  \end{eqnarray}
   and describes a pair of continuous curves enveloping a circle congruence.
\end{lem}
\noindent
{\bf Proof.} Consistency is shown exactly as in Lemma \ref{lem:dos}, the
limiting system is calculated directly from (\ref{eq:dP.})-(\ref{eq:dB.})
using (\ref{eq:Npositive}). One sees that the system of the lemma coincides
with (\ref{eq:ddd1})-(\ref{eq:ddd8}) for $M=2$. \BBox

\begin{prp} \label{prp:ribpairs}
{\bf (approximation of a pair of curves enveloping a
circle congruence).}
  Let $m=1$, $M=2$, and let there be given:
  \begin{itemize}
  \item a smooth curve $\XX:\PP_1\mapsto\rz^N$,
  \item a smooth function $A:\PP_1\mapsto\rz$,
  \item a point $\XX^+(0)\in\rz^N$.
  \end{itemize}
  Set $H_2(0)=|X^+(0)-X(0)|$, $\;\hat{v}_2(0)=H_2(0)^{-1}(X^+(0)-X(0))$.
  Assume that $\Psi\in\cH_\infty$ is suited for
  $(\hat{\XX}_1(0);\hat{v}_1(0))$, and set $\hat{v}_k(0)=A_{\Psi}(\ee_k)$
  for $k\neq 2$.
  Then, for some $r>0$:
  \begin{enumerate}
  \item There exists a unique curve $\XX^+:\PP_1\cap B_0(r)\mapsto\rz^N$
    through the point $\XX^+(0)$ such that the pair $(X,X^+)$ envelopes
    a circle congruence, and
    \begin{equation} \label{eq:getalpha}
     |\partial_1\XX^+|-|\partial_1\XX|=A\cdot|\XX^+-\XX|\,.
    \end{equation}
  \item Consider the family of pairs of discrete curves
  $(X^\eps,(X^\eps)^+)_{0<\eps<\eps_1}$ defined as the solutions
  to the Goursat problem \ref{gou:dos1} with
    the data
     \begin{eqnarray*}
      & \Psi^\eps=\Psi,\quad H^\eps_1=H_1\rest_{\PP_1^\eps},
     \quad B^\eps_{k1}=B_{k1}\rest_{\PP_1^\eps},
     \quad A^\eps=A\rest_{\PP_1^\eps}, & \\
      & H_2^\eps=H_2(0)\,,\quad
     B_{k2}^\eps=-2\langle\hat{v}_2(0),\hat{v}_k(0)\rangle. &
     \end{eqnarray*}
    where the functions $H_1$ and $B_{k1}$ are read off the smooth curve
    $X$, so that $X^\eps$ is the canonical
    $\eps$-discretization of $\XX$. These discrete curves
    $O(\eps)$-converge in $C^\infty$ to $(X,X^+)$.
  \end{enumerate}
\end{prp}
\noindent
{\bf Proof} follows from Lemma \ref{lem:dos1} and Theorem \ref{thm:approx}.
\BBox

\begin{thm}  \label{thm:dosapprox2}
{\bf (approximation of a Ribaucour pair).}
  Let, in addition to the data of Theorem \ref{thm:dosapprox}, there be
  given $m$ curves $\XX_i^+:\PP_i\mapsto\rz^N$ with a common intersection
  point $\XX^+=\XX^+_1(0)=\ldots=\XX^+_m(0)$, and such that each pair
  $(\XX_i,\XX_i^+)$ envelopes a circle congruence. In addition to the
  discrete surfaces ${\cal S}_{ij}^\eps$ from Theorem
  \ref{thm:dosapprox}, construct discrete curves $(X_i^\eps)^+$
  according to Proposition \ref{prp:ribpairs}, with the
  functions $A_i=(|\partial_iX_i^+|-|\partial_iX_i|)/|X_i^+-X_i|$.
  Then for some $r>0$:
  \begin{enumerate}
  \item
    There exists a unique Ribaucour pair of orthogonal systems
    $\xx,\xx^+:B_0(r)\mapsto\rz^N$ such that $\xx$ coincides with
    ${\cal S}_{ij}$ on $\PP_{ij}\cap B_0(r)$, and $\xx^+$
    coincides with $\XX_i^+$ on $\PP_i\cap B_0(r)$.
  \item
    There exists a unique family of Ribaucour pairs of discrete orthogonal
    systems $\{\xx^\eps:B^\eps(r)\mapsto\rz^N\}_{0<\eps<\eps_0}$ coinciding
    with ${\cal S}_{ij}^\eps$ on $(\PP_{ij}\times\{0\})\cap B(r)$ and
    coinciding with $(\XX_i^\eps)^+$ on $(\PP_i^\eps\times\{1\})\cap B(r)$.
    The family $x^\eps$ converges to the pairs $x,x^+$ with the rate
    $O(\eps)$ in $C^\infty$.
  \end{enumerate}
\end{thm}
\noindent
{\bf Proof} is similar to that of Theorem \ref{thm:dosapprox}, with
the only change in the conculding argument: for the limiting Jonas
pair $x,x^+$ of orthogonal nets we derive from eq. (\ref{circ}):
\[
  \hat{v}_M\hat{v}_i^++\hat{v}_i\hat{v}_M=0\quad\Rightarrow\quad
  \langle \hat{v}_i+\hat{v}_i^+,\hat{v}_M\rangle=0,
\]
which is the defining property of the Ribaucour pair. \BBox
\medskip

Letting not one but two or three directions of the orthogonal system
remain discrete in the continuous limit (so that $m'=2$ or $m'=3$), one
arrives at the following statement.
\begin{thm} \label{thm: rib permut}
{\bf (permutability of Ribaucour transformations).}
\begin{enumerate}
\item Given an $m$--dimensional orthogonal system
$x(\cdot,0,0):B_0(r)\mapsto\rz^N$ and its two Ribaucour transformations
$x(\cdot,1,0):B_0(r)\mapsto\rz^N$ and
$x(\cdot,0,1):B_0(r)\mapsto\rz^N$, there exists a
one--parameter family of orthogonal systems
$x(\cdot,1,1):B_0(r)\mapsto\rz^N$ that are Ribaucour transformations
of both $x(\cdot,1,0)$ and $x(\cdot,0,1)$. Corresponding points of
the four conjugate nets are concircular.
\item
Given three Ribaucour transformations
\[
x(\cdot,1,0,0),x(\cdot,0,1,0),x(\cdot,0,0,1):B_0(r)\mapsto\rz^N
\]
of a given $m$--dimensional orthogonal system
$x(\cdot,0,0,0):B_0(r)\mapsto\rz^N$, as well as three further
orthogonal systems
\[
x(\cdot,1,1,0),x(\cdot,0,1,1),x(\cdot,1,0,1):B_0(r)\mapsto\rz^N
\]
such that $x(\cdot,1,1,0)$ is a Ribaucour transformation of both
$x(\cdot,1,0,0)$ and $x(\cdot,0,1,0)$ etc., there exists
generically a unique orthogonal system
\[
x(\cdot,1,1,1):B_0(r)\mapsto\rz^N
\]
which is a Ribaucour transformation of all three $x(\cdot,1,1,0),
x(\cdot,0,1,1)$ and $x(\cdot,1,0,1)$.
\end{enumerate}
\end{thm}

\subsection{Example: elliptic coordinates}

\begin{figure}[t]
  \begin{center}
   \epsfig{file=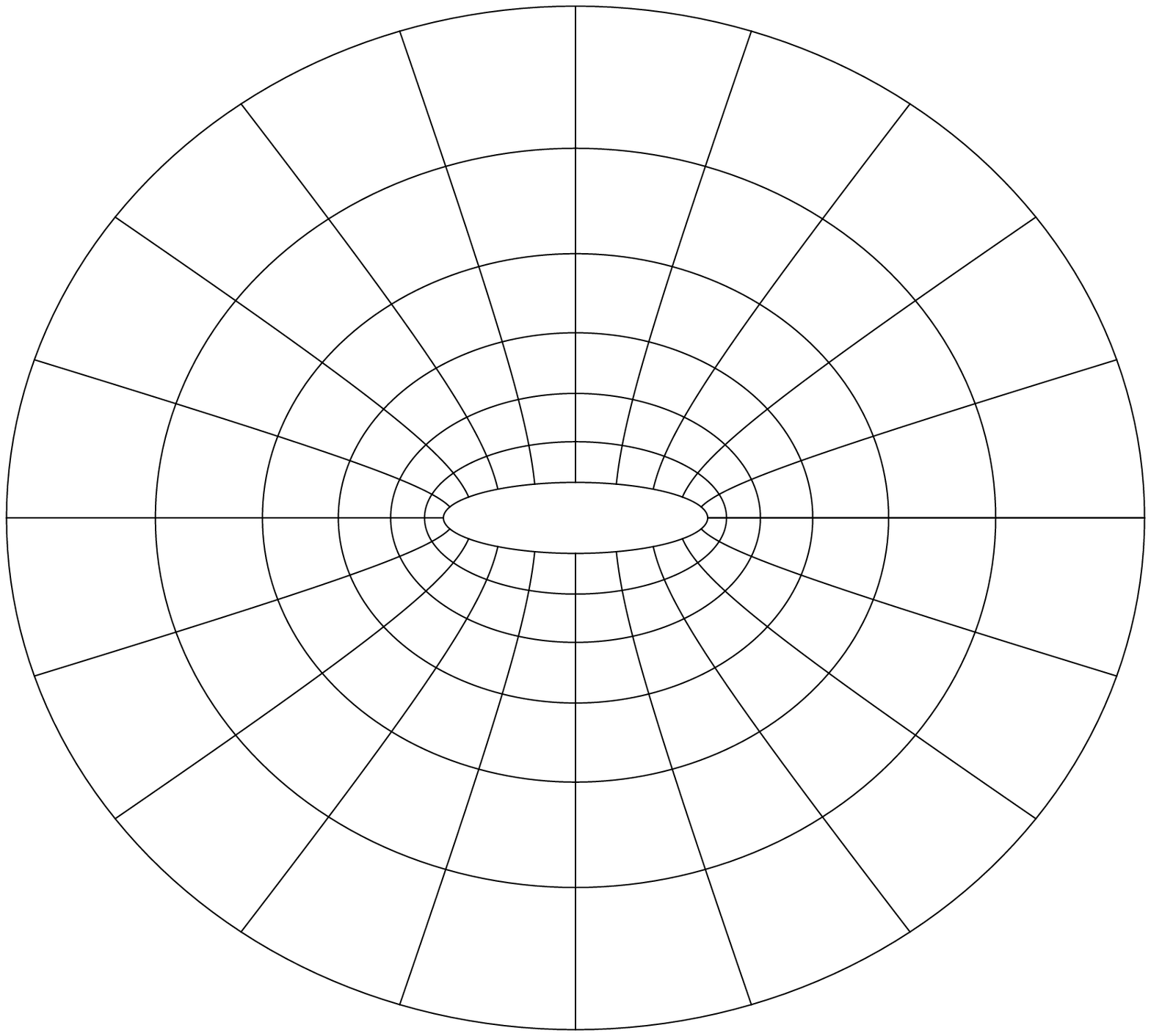,height=150pt}
   \epsfig{file=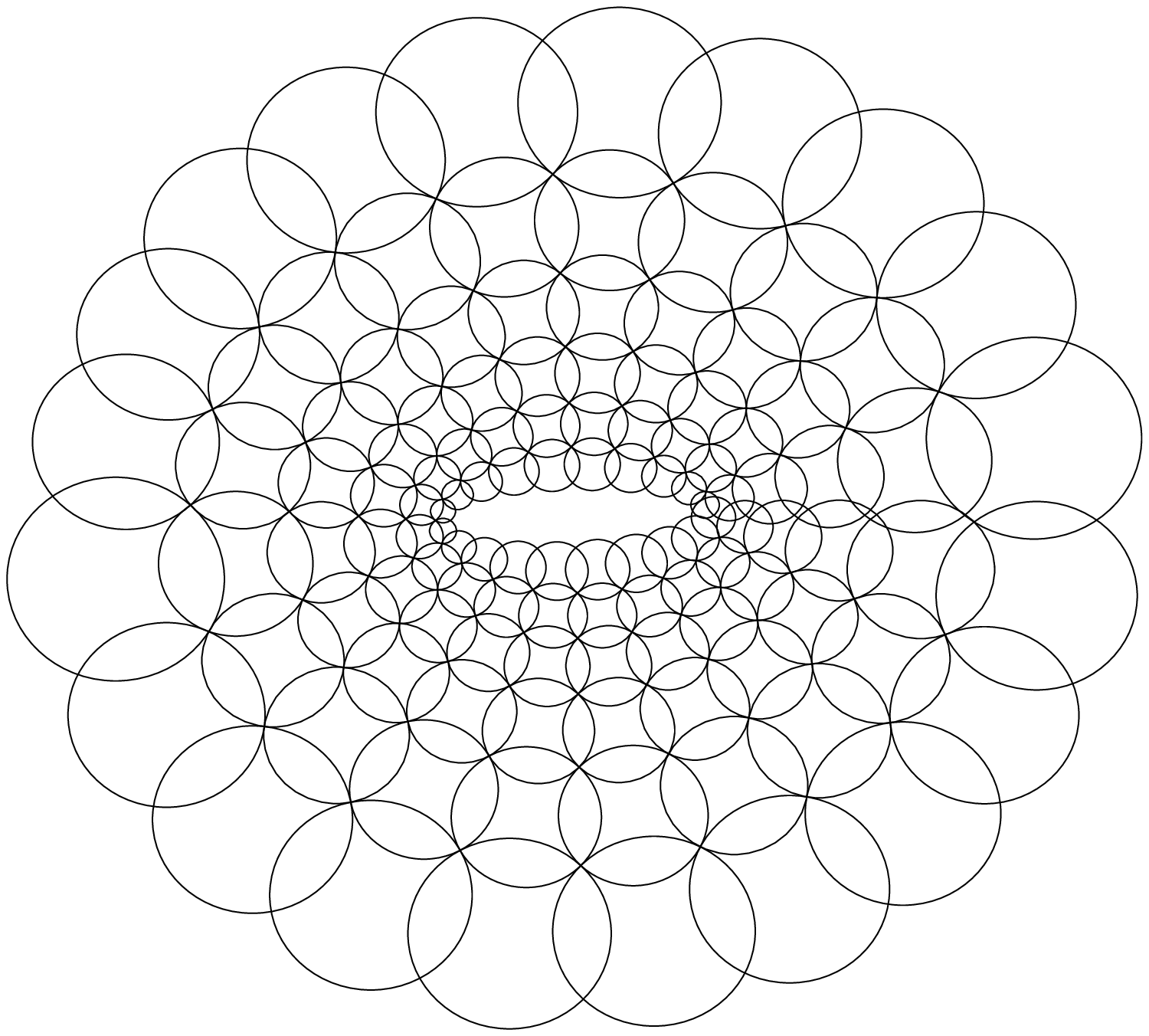,height=150pt}
   \vspace{5mm}\\
  \epsfig{file=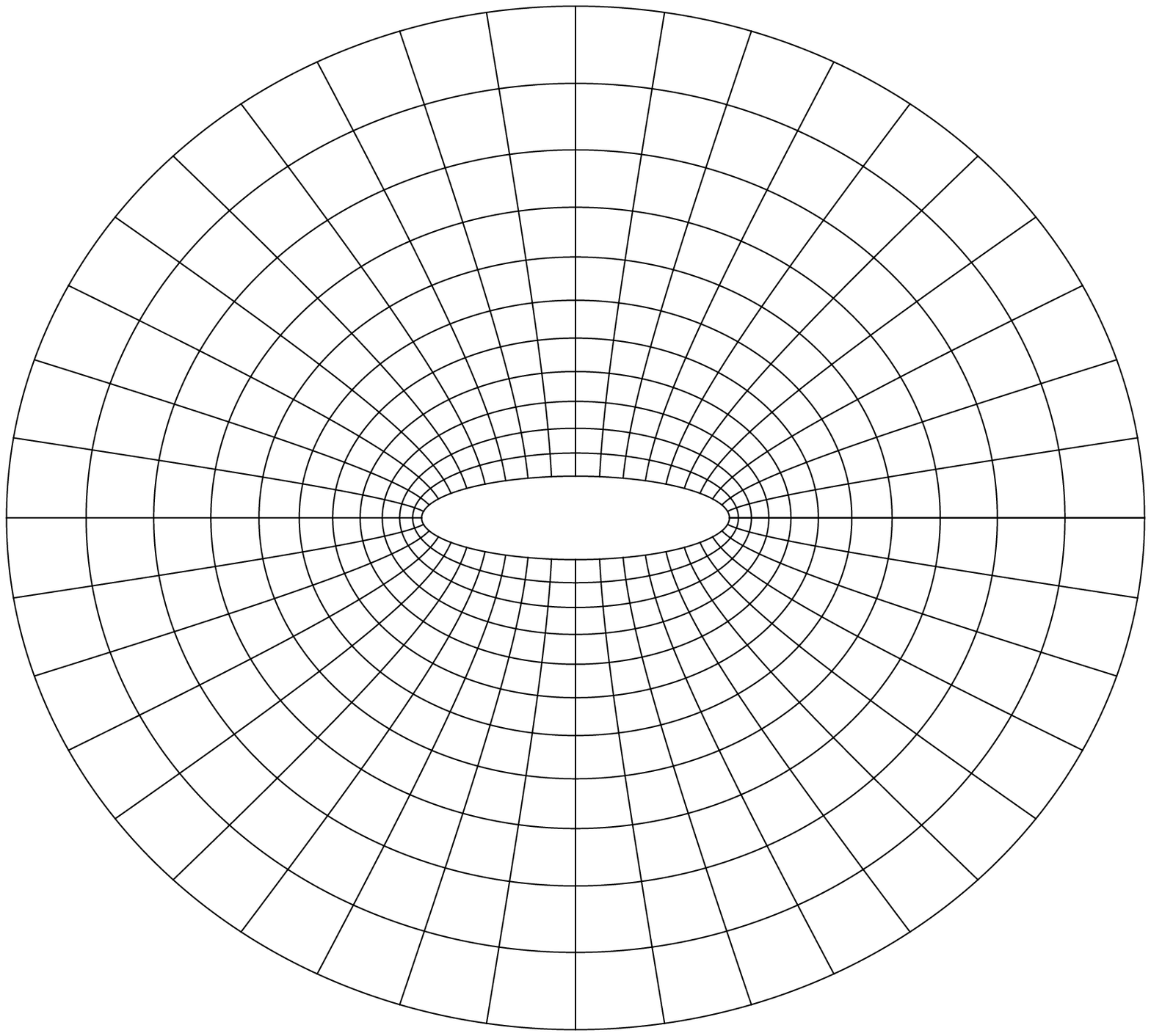,height=150pt}
  \epsfig{file=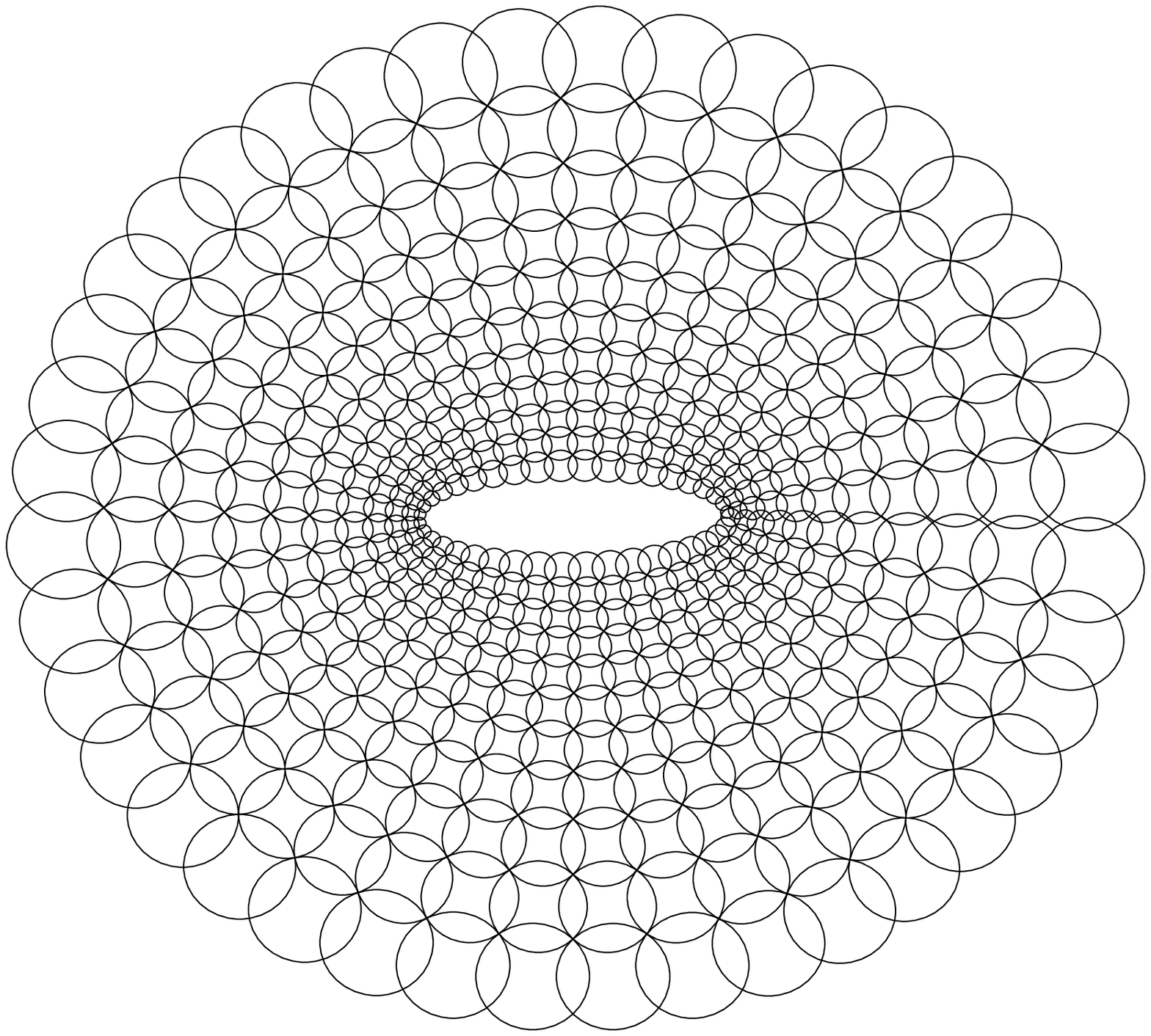,height=150pt}
  \end{center}
  \caption{Approximation of elliptic coordinates, $\eps=\pi/10$ and $\eps=\pi/20$.}
  \label{fig:few}
\end{figure}

The simplest nontrivial example, to which the above theory can be applied,
is the approximation of two-dimensional conformal maps by circular patterns.
Starting with a conformal map $F:\rz^2\rightarrow\rz^2$, i.e. $M=N=2$,
\begin{equation*}
  h=|\partial_1 F|=|\partial_2 F|\quad{\rm and}\quad
  \partial_1 F\cdot\partial_2 F=0,
\end{equation*}
one calculates the metric and rotation coefficients according to equations
(\ref{eq:DX}) and (\ref{eq:DH}), and the quantity $\gamma$ according to
(\ref{eq:DLI}):
\begin{equation*}
  h_1=h_2=h,\quad\beta_{12}=\frac{\partial_1 h}{h},\quad\beta_{21}=\frac{\partial_2 h}{h},\quad
  \gamma=\partial_1\beta_{12}=-\partial_2\beta_{21}
\end{equation*}
To construct a discrete approximation of $F$, solve the Goursat problem
\ref{gou:dos} for $\eps>0$ with the data
$H^\eps_i$, $B^\eps_{ij}$ and $\Gamma^\eps$ that are close to the
values of the respective functions $h_i$, $\beta_{ij}$ and $\gamma$
at the corresponding lattice sites. A good choice is, for example,
to prescribe
\begin{displaymath}
  \begin{array}{c}
    H^\eps_1(\xi_1)=h(\xi_1+\eps/2,0),\quad
    H^\eps_2(\xi_2)=h(0,\xi_2+\eps/2),\\
    B^\eps_{21}(\xi_1)=\beta_{21}(\xi_1+\eps/2,0),\quad
    B^\eps_{12}(\xi_2)=\beta_{12}(0,\xi_2+\eps/2),\\
    \Gamma^\eps(\xi_1,\xi_2)=\gamma(\xi_1+\eps/2,\xi_2+\eps/2),
  \end{array}
\end{displaymath}
for $\xi_1,\xi_2=i\eps$, $\;0\le i\le R$.
If $\psi^\eps:B^\eps(\eps R)\rightarrow\cH_\infty$ is the frame of
the respective solution to the system (\ref{eq:dP.})-(\ref{eq:dB.}), then
the function
\begin{equation*}
  F^\eps:B^\eps(R\eps)\mapsto\rz^2,\quad F^\eps(\xi_1,\xi_2)=
  \pi(\psi(\xi_1,\xi_2)^{-1}\ee_0\psi(\xi_1,\xi_2))
\end{equation*}
is a two-dimensional discrete orthogonal system, i.e., the points
$F^\eps(\xi_1,\xi_2)$, $F^\eps(\xi_1+\eps,\xi_2)$, $F^\eps(\xi_1,\xi_2+\eps)$
and $F^\eps(\xi_1+\eps,\xi_2+\eps)$ lie on a common circle
${\rm C}^\eps(\xi_1,\xi_2)$ in $\rz^2$,
and $F^\eps(\xi_1,\xi_2)=F(\xi_1,\xi_2)+O(\eps)$.

This is illustrated with the planar elliptic coordinate system:
\begin{equation*}
  F(\xi_1,\xi_2)=(\cosh(\xi_1)\cos(\xi_2),
  \sinh(\xi_1)\sin(\xi_2)),
\end{equation*}
whose coordinate lines are ellipses and hyperbolas.
One finds
\begin{eqnarray*}
 & h=\left(\sinh^2(\xi_1)+\sin^2(\xi_2)\right)^{1/2},& \\
 & \beta_{12}=\sinh(2\xi_1)/(2h^2),\quad
  \beta_{21}=\sin(2\xi_2)/(2h^2), &\\
 & \gamma=\left(1-\cosh(2\xi_1)\cos(2\xi_2)\right)/(4h^4). &
\end{eqnarray*}
Results are displayed in Fig.
\ref{fig:few}.
Their left sides show the original coordinate lines of $F$,
and on the right sides the circles ${\rm C}^\eps(\xi_1,\xi_2)$ are drawn;
each intersection point of two coordinate lines on the left
corresponds to an intersection point of four circles on the right.
There are small defects (the circles do not close up) on the very right
of the pictures since, in contrast to $F$, the discrete maps $F^\eps$
are not periodic with respect to $\xi_2$.


\end{document}